 \documentclass[12pt]{article}
 
 \usepackage{amssymb}		
 \usepackage{amsmath}		
 \usepackage{latexsym}		
 \usepackage{epsfig}         	
 \usepackage{subfigure}
 \usepackage{pstricks,pst-node,pst-text,pst-tree}

 \newtheorem{definition}{Definition}[section]
 
 \newtheorem{corollary}[definition]{Corollary}
 \newtheorem{example}[definition]{Example}

 \newtheorem{lemma}[definition]{Lemma}
 \newtheorem{proposition}[definition]{Proposition}
 \newtheorem{remark}[definition]{Remark}
 \newtheorem{theorem}[definition]{Theorem}

 \numberwithin{equation}{section}

 \newcommand{\proof}{\noindent\emph{Proof.\ }}
 \newcommand{\qed}{\hspace*{\fill}\(\Box\) \medskip}
 \newcommand{\one}{\ensuremath{(\mathrm{i})}}
 \newcommand{\two}{\ensuremath{(\mathrm{ii})}}
 \newcommand{\three}{\ensuremath{(\mathrm{iii})}}

 \newcommand{\C}{\ensuremath{\mathbb{C}}}
 \newcommand{\Div}{\ensuremath{\operatorname{Div}}}
 
 \newcommand{\Hom}{\ensuremath{\mbox{Hom}}}
 
 \newcommand{\Pic}{\ensuremath{\operatorname{Pic}}}
 \newcommand{\Q}{\ensuremath{\mathbb{Q}}}
 \newcommand{\R}{\ensuremath{\mathbb{R}}}
 \newcommand{\Region}{\ensuremath{\operatorname{Conv}}}
 \newcommand{\SL}{\ensuremath{\operatorname{SL}}}
 \newcommand{\Spec}{\ensuremath{\operatorname{Spec}}}
 \newcommand{\Z}{\ensuremath{\mathbb{Z}}}

 \newcommand{\ahilb}{\ensuremath{A\operatorname{-Hilb}(\mathbb{C}^{3})}}

 \newcommand{\dimn}{\operatorname{dim}}
 \newcommand{\divisor}{\operatorname{div}}
 \newcommand{\dP}{\ensuremath{\operatorname{dP_{6}}}}
 
 \newcommand{\owe}{\ensuremath{\mathcal{O}}}
 \newcommand{\rank}{\ensuremath{\operatorname{rank}}}
 \newcommand{\st}{\ensuremath{\operatorname{\bigm{|}}}}

 \newcommand{\wt}{\ensuremath{\operatorname{wt}}}

 \newcommand{\Ahilb}[1]{\ensuremath{A\operatorname{-Hilb}(\mathbb{C}^{#1})}}
 \newcommand{\Ghilb}[1]{\ensuremath{G\operatorname{-Hilb}(\mathbb{C}^{#1})}}

\makeatletter
 \newlength{\typesize}
\makeatother

\newlength{\vvoff}
\newlength{\hhoff}

\newcommand{\locateoffcenter}[1]{%
\addtolength{\vvoff}{-0.25\typesize}%
\raisebox{\vvoff}{\hspace{\hhoff}\makebox(0,0){\smash{#1}}}
}
\newcommand{\object}[1]{%
\setlength{\vvoff}{0pt}%
\setlength{\hhoff}{0pt}%
\locateoffcenter{#1}
}

\newcommand{\swlabel}[1]{%
\setlength{\vvoff}{-0.5\typesize}%
\setlength{\hhoff}{0.75\typesize}%
\locateoffcenter{#1}
}
\newcommand{\nwlabel}[1]{%
\setlength{\vvoff}{-0.5\typesize}%
\setlength{\hhoff}{-0.75\typesize}%
\locateoffcenter{#1}
}
\newcommand{\selabel}[1]{%
\setlength{\vvoff}{0.5\typesize}%
\setlength{\hhoff}{0.75\typesize}%
\locateoffcenter{#1}
}
\newcommand{\nelabel}[1]{%
\setlength{\vvoff}{0.5\typesize}%
\setlength{\hhoff}{-0.75\typesize}%
\locateoffcenter{#1}
}

 \title{An explicit construction of the McKay correspondence for \protect\(A\protect\)-Hilb \protect\(\C^3\protect\)}
 \author{Alastair Craw}
 \date{}

\begin{document}
 
 \maketitle

 \begin{abstract} 
 For a finite Abelian subgroup \(A \subset \SL(3,\C)\),  let \(Y = \ahilb\) denote the scheme parametrising \(A\)-clusters in \(\C^{3}\).  Ito and Nakajima proved that the tautological line bundles (indexed by the irreducible representations of \(A\)) form a basis of the \(K\)-theory of \(Y\).  We establish the relations between these bundles in the Picard group of \(Y\) and hence,  following a recipe introduced by Reid,  construct an explicit basis of the integral cohomology of \(Y\) in one-to-one correspondence with the irreducible representations of \(A\).  

 MSC 2000: primary 14E15; secondary 14F05, 14J30, 19L64.
 \end{abstract}

 \section{Introduction}
 \label{sec:intro}

 Let \(G\subset \SL(n,\C)\) be a finite subgroup.   A \emph{\(G\)-cluster} is a \(G\)-invariant zero-dimensional subscheme  \(Z \subset \C^{n}\) with global sections \(H^{0}(Z,\mathcal{O}_{Z})\) isomorphic as a \(\C[G]\)-module to the regular representation of \(G\).  Write \(\Ghilb{n}\) for the moduli space of \(G\)-clusters.  Ito and Nakamura~\cite{Ito:hs} proved that \(\Ghilb{2}\) is the unique minimal (or \emph{crepant}) resolution \(Y\) of \(\C^{2}/G\).  Nakamura~\cite{Nakamura:ago} conjectured that \(\Ghilb{3}\) is a crepant resolution of the quotient $\C^3/G$ and proved this for a finite Abelian subgroup $A\subset \SL(3,\C)$ by introducing an algorithm that calculates \ahilb.  Nakamura's conjecture was subsequently proved by Bridgeland,  King and Reid~\cite{Bridgeland:mim}.

 The search for crepant resolutions of \(\C^{n}/G\) was motivated in part by the McKay correspondence.  For a finite subgroup \(G \subset \SL(2,\C)\),  McKay~\cite{McKay:gsfg} established a one-to-one correspondence between the nontrivial irreducible representations of \(G\) and the exceptional prime divisors of the crepant resolution \(Y\) of \(\C^{2}/G\).  Gonzalez-Sprinberg and Verdier~\cite{Gonzalez-Sprinberg:cgm} subsequently provided a geometric explanation by associating a vector bundle \(\mathcal{R}_{k}\) on \(Y\) to each  irreducible representation \(\rho_{k}\) of \(G\).  Case by case analysis of the finite subgroups \(G \subset \SL(2,\C)\) revealed that the classes \(c_{1}(\mathcal{R}_{k})\) (for nontrivial \(\rho_{k}\)) form a basis of \(H^{2}(Y,\Z)\) dual to the exceptional divisor classes.  This leads to a one-to-one correspondence
 \begin{equation}
 \label{eqn:bij}
 \Big{\{}\mbox{irreducible representations of\ }G\Big{\}}\; \longleftrightarrow\; \mbox{basis of\ } H^*(Y,\mathbb{Z}).
 \end{equation}

 Following Nakamura's conjecture that \(\Ghilb{3}\) is a crepant resolution of \(\C^{3}/G\),  Reid~\cite{Reid:mc} conjectured that the tautological bundles \(\mathcal{R}_{k}\) on \(Y\) (see \S\ref{sec:tlb} for the definition) form a \(\Z\)-basis of the \(K\)-theory of \(\Ghilb{3}\),  and that a certain cookery with the Chern classes \(c_{i}(\mathcal{R}_{k})\) gives a \(\Z\)-basis of the integral cohomology of \(\Ghilb{3}\) satisfying (\ref{eqn:bij}).  To support the conjecture,  Reid calculated \(Y \! = \ahilb\) for several examples of finite Abelian subgroups \(A\subset \SL(3,\C)\) and decorated the toric fan of \(Y\) with characters of \(A\) in a manner that encoded the relations in \(\Pic(Y)\) between the line bundles \(\mathcal{R}_{k}\).  Every such relation in \(\Pic(Y)\) led to the construction of a virtual bundle \(\mathcal{V}_{m}\) indexed by an irreducible representation \(\rho_{m}\) of \(A\).  In each of Reid's examples,   the second Chern classes \(c_{2}(\mathcal{V}_{m})\) base \(H^{4}(Y,\Z)\) and,  moreover, the first Chern classes \(c_{1}(\mathcal{R}_{k})\) indexed by the remaining nontrivial irreducible representations \(\rho_{k}\) form a \(\Z\)-basis for \(H^{2}(Y,\Z)\).  Thus,  the McKay correspondence (\ref{eqn:bij}) holds for each of Reid's examples.

  Ito and Nakajima~\cite{Ito:mc3} proved the first part of Reid's conjecture for a finite Abelian subgroup \(A\subset \SL(3,\C)\), i.e., that the line bundles \(\mathcal{R}_{k}\) form a \(\Z\)-basis of the \(K\)-theory of \(Y\).  Applying the Chern character provides a basis of \(H^{*}(Y,\Q)\) in one-to-one correspondence with the irreducible representations of \(A\),   a rational version of (\ref{eqn:bij}).  In this paper we establish the integral version of (\ref{eqn:bij}) for a finite Abelian subgroup \(A\subset \SL(3,\C)\):

 \begin{theorem}
 \label{thm:McKay}
 The McKay correspondence bijection (\ref{eqn:bij}) holds (replace \(G\) by \(A\)) for all finite Abelian subgroups \(A\subset \SL(3,\C)\).
 \end{theorem}

 The first step is to show that the recipe introduced by Reid~\cite{Reid:mc} which decorates the lines and vertices in (a cross-section of) the toric fan \(\Sigma\) of \ahilb\ with characters of the group \(A\) can be carried out for any finite Abelian subgroup \(A\subset \SL(3,\C)\).  In addition, we prove that every character of \(A\) marks either a line in \(\Sigma\) (possibly passing through several vertices) or a unique vertex.  See \S\ref{sec:df} for examples.  

 The decoration of \(\Sigma\) with characters enables us to calculate the relations between the line bundles \(\mathcal{R}_{k}\) in \(\Pic(Y)\).  For each interior vertex \(v\) in \(\Sigma\),  we derive a relation between those bundles \(\mathcal{R}_{k}\) indexed by the irreducible representations \(\rho_{k}\) whose characters mark the vertex \(v\) and the lines meeting at \(v\) (see Theorem~\ref{thm:bundles} for the precise statement and a list of the explicit relations).  A weak version of the McKay correspondence,  namely the equality of the Euler number of \(Y\) and the order of the group \(A\),  shows that our list exhausts all nontrivial relations between tautological line bundles.   These relations are of independent interest,  see Craw and Ishii~\cite[\S7]{Craw:fog}.

 Once the relations in \(\Pic(Y)\) have been derived,  we implement Reid's construction of virtual bundles \(\mathcal{V}_{m}\) on \(Y\) having trivial rank and trivial first Chern class.  A proof based on a case-by-case analysis of the relations in \(\Pic(Y)\) shows that the second Chern classes \(c_{2}(\mathcal{V}_{m})\) of these virtual bundles  form a basis of \(H^{4}(Y,\Z)\) dual to the basis \([S]\in H_{4}(Y,\Z)\) of the compact exceptional surfaces \(S\) of the resolution \(\varphi\colon Y \to \C^{3}/A\).   To construct the \(\Z\)-basis of \(H^{2}(Y,\Z)\),  we start with the spanning set of all first Chern classes \(c_{1}(\mathcal{R}_{k})\) of tautological bundles.  Removing those classes \(c_{1}(\mathcal{R}_{m})\) indexed by the irreducible representations \(\rho_{m}\) corresponding to the virtual bundles \(\mathcal{V}_{m}\) leaves a \(\Z\)-basis of \(H^{2}(Y,\Z)\).  Theorem~\ref{thm:McKay} then follows since the trivial bundle \(\mathcal{R}_{\rho_{0}} = \owe_{Y}\) generates \(H^{0}(Y,\Z)\).

 The work of Bridgeland,  King and Reid~\cite{Bridgeland:mim} implies that the bundles \(\mathcal{R}_{k}\) form a \(\Z\)-basis of the \(K\)-theory of \(Y\) for any finite subgroup \(G\subset \SL(3,\C)\),  but the cookery with the Chern classes leading to a \(\Z\)-basis of \(H^{*}(Y,\Z)\),  is still an open problem for non-Abelian subgroups \(G\subset \SL(3,\C)\).
 
 \medskip

 \noindent \textbf{Acknowledgements\ } I wish to thank M.~Reid for encouraging me to investigate the McKay correspondence in light of \cite{Craw:htc}.  I am grateful to M.~Brion and G.~Gonzalez-Sprinberg for pointing out a gap in an earlier version of this article.  Thanks also to Y.~Ito,  A.~King, R.~Leng and B.~Szendr\H{o}i.

 \section{How to calculate \ahilb}
 \label{sec:htc}
 Let \(A \subset \SL(3,\C)\) be a finite Abelian subgroup of order \(r = \vert A\vert\),  and fix a primitive \(r\)th root of unity \(\varepsilon\).  Choose coordinates \(x,y,z\) on \(\C^{3}\) to diagonalise the action of \(A\),  write \(L \cong \Z^3\) for the lattice of Laurent monomials in \(x,y,z\) and \(L^{\vee}\) for the dual lattice with basis \(e_{1},e_{2},e_{3}\).  To each group element \(a = \mbox{diag}(\varepsilon^{\alpha_{1}},\varepsilon^{\alpha_{2}},\varepsilon^{\alpha_{3}})\) with \(0\leq \alpha_{j} < r\),  we associate the vector \(v_{a} = \frac{1}{r}(\alpha_{1},\alpha_{2},\alpha_{3})\).  Write \(N := L^{\vee} + \sum_{a\in A} \; \Z\cdot v_{a}\) and \(M := \Hom(N,\Z)\) for the dual lattice of \(A\)-invariant Laurent monomials.  

 The toric variety \(U_{\sigma} := \Spec\ \C[\sigma^{\vee}\cap M]\) defined by the positive octant \(\sigma = \sum \R_{\geq 0} e_{i}\) in \(N_{\R} := N\otimes_{\Z}\R\) is the quotient \(\C^3/A\).  The \emph{junior simplex} \(\Delta \subset N_{\R}\) is the triangle with vertices \(e_{1}, e_{2}, e_{3}\),  containing the lattice points \(\frac{1}{r}(\alpha_{1}, \alpha_{2}, \alpha_{3})\) with \(\frac{1}{r}(\alpha_{1} + \alpha_{2} + \alpha_{3}) = 1\).  Crepant toric resolutions of \(\C^{3}/A\) are determined by triangulations of the junior simplex \(\Delta\) into basic triangles (we identify a triangulation of \(\Delta\) with the fan determining the triangulation).  Nakamura~\cite{Nakamura:ago} exhibited one such triangulation \(\Sigma\) determining the toric variety \(X_{\Sigma} = \Ahilb{3}\) parametrising \(A\)-clusters.  An \emph{\(A\)-cluster} is an \(A\)-invariant,  zero-dimensional subscheme \(Z \subset \C^{3}\) for which \(H^{0}(Z,\owe_{Z})\) is isomorphic as a \(\C[A]\)-module to the regular representation of \(A\).  Craw and Reid~\cite{Craw:htc} calculated \(\Sigma\) by the following three-step procedure:
 \begin{enumerate}
 \item Draw lines \(L_{i,0},\dots ,L_{i,m_{i}+1}\) emanating from the corners \(e_{i}\) of \(\Delta\) to the points forming the convex hull of lattice points in \(\Delta \smallsetminus e_{i}\) (the lines \(L_{i,0}\) and \(L_{i,m_{i}+1}\) extend along two sides of \(\Delta\)).  For \(j = 1,\dots ,m_{i}\) the integer \(a_{i,j}\geq 2\) determined by the Jung--Hirzebruch continued fraction rule \(L_{i,j-1} + L_{i,j+1} = a_{i,j}\cdot L_{i,j}\) is called the \emph{strength} of \(L_{i,j}\).

 \item Extend the lines \(L_{i,1},\dots ,L_{i,m_{i}}\) until they are `defeated' by lines \(L_{k,l}\) from \(e_{k}\) (\(i\neq k\)) according to the following rule:  when lines meet at a point,  the line with greatest strength extends with strength reduced by 1 for every rival it defeats; lines meeting with equal strength all die.  This results in a partition of \(\Delta\) into \emph{regular triangles} of side \(r\), i.e.,  lattice triangles with \(r+1\) lattice points on each edge.  
 \item Draw \(r-1\) lines parallel to the sides of each regular triangle of side \(r\) to produce its \emph{regular tesselation} into \(r^{2}\) basic triangles (see Craw and Reid~\cite[Figure~2(a)]{Craw:htc}).  The resulting basic triangulation is \(\Sigma\).
 \end{enumerate}
 
\begin{figure}[!ht]
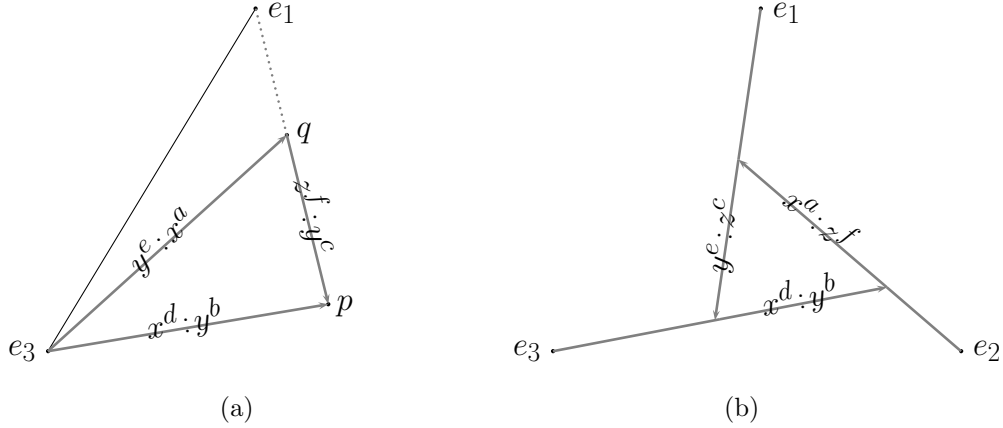

     \centering
     \mbox{\subfigure[]{\input{regcorner.pspic}\label{fig:corner_triangle}} \hspace{1.3cm}
           \subfigure[]{\input{regchamp.pspic}\label{fig:champion_triangle}}}
     \caption{(a) corner triangle $e_{3}pq$; (b) meeting of champions} 
     \label{fig:reg_triangles} 
 \end{figure}

 To prove that \(X_{\Sigma} = \ahilb\),  we pass to the dual \(M = \Hom(N,\Z)\) and calculate coordinates on the affine pieces \(U_{\tau} = \Spec\ \C[\tau^{\vee}\cap M]\) covering \(X_{\Sigma}\).  In doing so we prove that every regular triangle of side \(r\) is either a \emph{corner triangle} or a (unique) \emph{meeting of champions} as shown in Figure~\ref{fig:reg_triangles} (permute \(x,y,z\) if necessary). 
 The indices of the \(A\)-invariant ratios cutting out the sides of the regular triangles in Figure~\ref{fig:reg_triangles} satisfy
\begin{eqnarray}
     d - a = e - b - c = f = r \quad \mbox{in Case a,} \label{eqn:prop3.1a}  \\
     d - a = e - b = f - c = r \quad \mbox{in Case b.} \label{eqn:prop3.1b}
\end{eqnarray}  
Moreover,  the lines of the regular tesselations of the regular triangles of Figure~\ref{fig:reg_triangles} are cut out by the \(A\)-invariant ratios of monomials
\begin{eqnarray}
     x^{d-i}\!:\!y^{b+i}z^{i};\quad y^{e-j}\!:\!z^{j}x^{a+j};\quad z^{f-k}\!:\!x^{k}y^{c+k}\quad \mbox{in Case a,} \label{eqn:prop3.2a}\\
     x^{d-i}\!:\!y^{b+i}z^{i},\quad y^{e-j}\!:\!z^{c+j}x^{j},\quad z^{f-k}\!:\!x^{a+k}y^{k}\quad \mbox{in Case b,}\label{eqn:prop3.2b}
\end{eqnarray}
for \(i,j,k = 0,\dots ,r-1\).  The edges of a basic triangle \(\tau\in \Sigma\) are cut out by the ratios from (\ref{eqn:prop3.2a}) or (\ref{eqn:prop3.2b}) if \(i,j,k\) satisfy either \(i+j+k = r-1\),  in which case \(\tau\) is said to be an \emph{up} triangle,  or \(i+j+k=r+1\),  in which case \(\tau\) is \emph{down}.  An up triangle \(\tau\) defines the affine subvariety \(U_{\tau} \subset X_{\Sigma}\) isomorphic to \(\C^{3} = \Spec\ \C[\xi,\eta,\zeta]\) with coordinates
  \begin{align}
 & \xi = x^{d-i}/y^{b+i}z^{i},\quad \eta = y^{e-j}/z^{j}x^{a+j}, \quad \zeta = z^{f-k}/x^{k}y^{c+k} & \mbox{in Case a,} \label{eqn:coordsupa}\\
 & \xi = x^{d-i}/y^{b+i}z^{i},\quad \eta = y^{e-j}/z^{c+j}x^{j}, \quad \zeta = z^{f-k}/x^{a+k}y^{k} & \mbox{in Case b.} \label{eqn:coordsupb}
 \end{align}
Similarly,  a down triangle \(\tau\) defines the affine subvariety \(U_{\tau} \subset X_{\Sigma}\) which is isomorphic to \(\C^{3} = \Spec\ \C[\lambda,\mu,\nu]\) with coordinates
  \begin{align}
 & \lambda = y^{b+i}z^{i}/x^{d-i},\quad \mu = z^{j}x^{a+j}/y^{e-j}, \quad \nu = x^{k}y^{c+k}/z^{f-k} & \mbox{in Case a,} \label{eqn:coordsdowna}\\
 & \lambda = y^{b+i}z^{i}/x^{d-i},\quad \mu = z^{c+j}x^{j}/y^{e-j}, \quad \nu = x^{a+k}y^{k}/z^{f-k} & \mbox{in Case b.} \label{eqn:coordsdownb}
 \end{align}
 To complete the proof that \(X_{\Sigma} = \ahilb\),  it remains to compare the coordinates on the subsets \(U_{\tau}\) covering \(X_{\Sigma}\) with the explicit coordinates on a cover of \ahilb\ calculated by Nakamura~\cite{Nakamura:ago}. See \cite[\S4--5]{Craw:htc} for more details.

 \begin{remark}
 \label{remark:knock-out}
 \emph{The knock-out rule (Step 2) in the calculation of \(\Sigma\) can be given in terms of monomials. Indeed,  suppose a line \(L_{1,j}\) meets a line \(L_{3,k}\) at an interior point of \(\Delta\).  The lines are cut out by \(y^{c}\! : \! z^{f}\) and \(x^{a}\! : \! y^{e}\) respectively.  The knock-out rule can be stated as:  \emph{a line extends if and only if its ratio contains the strictly smaller exponent of the common monomial \(y\)}.  For example,  in Figure~\ref{fig:corner_triangle} above,   lines cut out by ratios \(y^{c}\! : \! z^{f}\) and \(x^{a}\! : \! y^{e}\) meet at a point an interior point of \(\Delta\).  The former line extends,  so \(c < e\).}
 \end{remark}

 \begin{example}
 \label{ex:11.1.2.8}
 \emph{Consider the cyclic quotient singularity of type \(\frac{1}{11}(1,2,8)\).  In Figure~\ref{fig:step1} we illustrate the result of Step~1 where,  for example,  the strengths of the lines from \(e_{3}\) come from the continued fraction \(\frac{11}{2} = 6 - \frac{1}{2}\) of the surface singularity \(\C^{2}_{(z=0)}/A = \frac{1}{11}(1,2)\).  
 \begin{figure}[!ht]
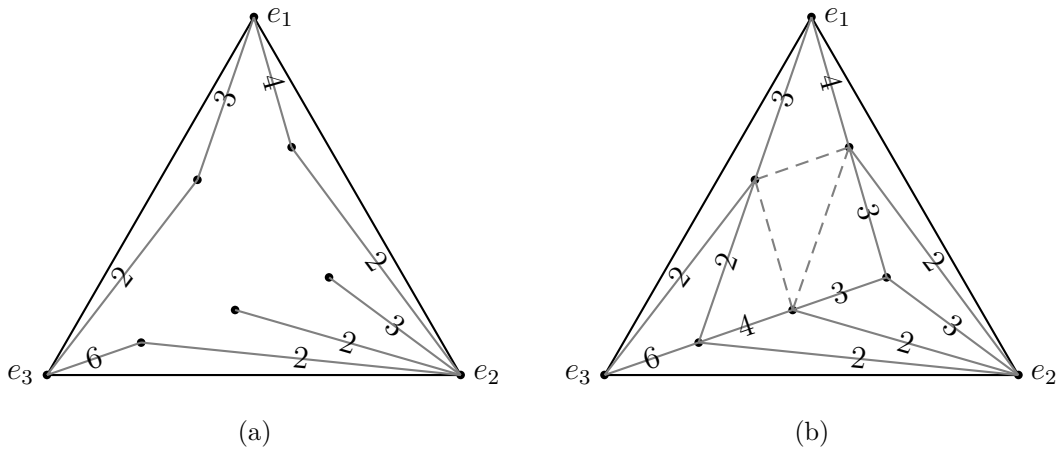

 \centering
 \mbox{\subfigure[]{\input{step1.pspic}\label{fig:step1}} \hspace{1.5cm} \subfigure[]{\input{step2.pspic}\label{fig:step2}}}
 \caption{(a) Step 1; (b) Step 2 (solid lines) and Step 3 (dotted lines)}
 \label{fig:steps}
 \end{figure}  
 The solid lines in Figure~\ref{fig:step2} show the result of Step~2.  The line from \(e_{1}\) with strength 3 intersects the line from \(e_{3}\) with strength 2,  so the line from \(e_{1}\) extends with strength 2 while the line from \(e_{3}\) terminates.  The resulting partition of \(\Delta\) contains only one regular triangle of side \(r > 1\).  To perform Step~3,  tesselate this triangle,  i.e.,  add the dotted lines to Figure~\ref{fig:step2},   giving \(\Sigma\).  The \(A\)-invariant ratios that cut out the lines in \(\Sigma\) are shown in Figure~\ref{fig:ratios} from \S\ref{sec:df}.}
 \end{example}

 \section{Reid's recipe for decorating \protect\(\Sigma\protect\)}
 \label{sec:df}
 Reid~\cite{Reid:mc} calculated \(X_{\Sigma} = \ahilb\) for several examples and,  in each case,  marked the lines and vertices in \(\Sigma\) with characters of \(A\).  We now prove that this can be done for any finite Abelian subgroup \(A\subset \SL(3,\C)\).

 Lines in \(\Sigma\) are cut out by the \(A\)-invariant ratios listed in (\ref{eqn:prop3.2a}) and (\ref{eqn:prop3.2b}).   The monomials in each ratio lie in the same character space of the \(A\)-action,  and we mark the line with the common character.  As for the vertices,  Reid introduced a recipe to associate one or two characters to a vertex \(v\),  depending (primarily) on the valency of \(v\).   In light of \cite[\S1.3]{Craw:htc},  the valency is either 3, 4, 5 or 6.  We now implement Reid's recipe case by case.

 \medskip

 \noindent\textsc{Case 1:} A vertex \(v\) of valency 3 defines an exceptional \(\mathbb{P}^{2}\).
 \begin{lemma}
 \label{lemma:val3}
 A single character \(\chi_{k}\) marks all three lines meeting at \(v\).  Mark the vertex \(v\) with the character \(\chi_{m} \! := \chi_{k}\otimes \chi_{k}\).
 \end{lemma}
 \proof
 A vertex \(v\) of valency 3 occurs only when a line \(L_{1,\alpha}\) emanating from \(e_{1}\) meets lines \(L_{2,\beta}\) and \(L_{3,\gamma}\) from \(e_{2}\) and \(e_{3}\) respectively.  The ratio cutting out \(L_{1,\alpha}\) is of the form \(y^{b}\! : \! z^{f}\).  The lines \(L_{1,\alpha}\) and \(L_{3,\gamma}\) defeat each other at \(v\) so,  by Remark~\ref{remark:knock-out},  the ratio cutting out \(L_{3,\gamma}\) is of the form \(x^{d}\! : \! y^{b}\).  Similarly,  \(L_{2,\beta}\) is cut out by \(z^{f}\! : \! x^{d}\).  In particular,  all three lines are marked with the common character of \(x^{d},  y^{b}, z^{f}\). \qed 

\begin{figure}[!ht]
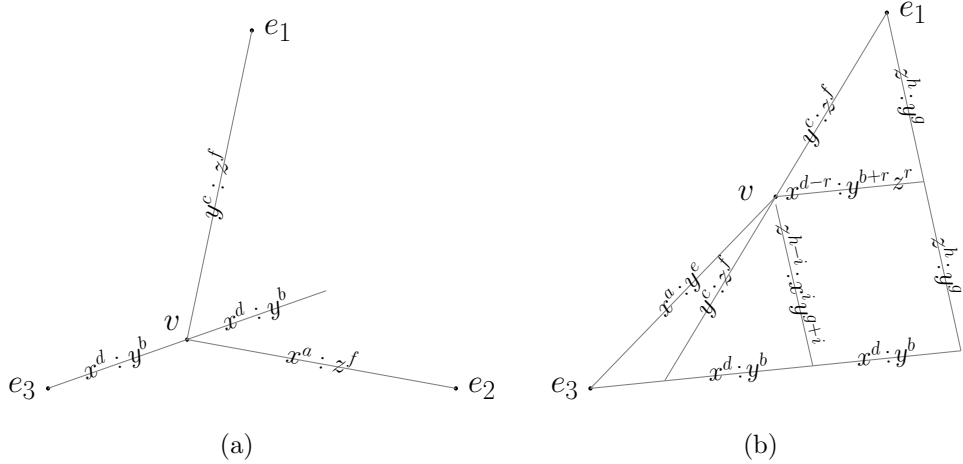

 \centering
 \mbox{\subfigure[]{\input{Fval4.pspic}\label{fig:Fval4}} \hspace{1.8cm} \subfigure[]{\input{Fval5.pspic}\label{fig:Fval5}}}
 \caption{Ratios on lines meeting at vertices of valency 4 and 5}
 \end{figure}

 \noindent\textsc{Case 2:} A vertex \(v\) of valency 4 defines an exceptional scroll \(\mathbb{F}_{r}\).
 \begin{lemma}
 \label{lemma:val4}
 There are distinct characters \(\chi_{k}\) and \(\chi_{l}\) which each mark a pair of lines meeting at \(v\).  Mark the vertex \(v\) with the character \(\chi_{m} \! := \chi_{k}\otimes \chi_{l}\).
 \end{lemma}
 \proof
 A vertex \(v\) of valency 4 occurs only when a line \(L_{\alpha,\beta}\) from \(e_{\alpha}\) defeats lines emanating from both of the other corners of \(\Delta\).  By permuting \(x,y,z\) if necessary we assume that \(\alpha =3\) (see Figure~\ref{fig:Fval4}).  Let \(\chi_{k}\) denote the common character space of the monomials \(x^{d}\) and \(y^{b}\) in the ratio marking \(L_{3,\beta}\).  If there are no vertices on \(L_{3,\beta}\) between \(e_{3}\) and \(v\) then \(z\) is one of the monomials in the ratios marking the defeated lines from \(e_{1}\) and \(e_{2}\).  More generally,  it follows from the calculation of \(\Sigma\) that if \((f\!-\!1)\) vertices lie between \(e_{3}\) and \(v\) then \(z^{f}\) occurs in both \(A\)-invariant ratios on the defeated lines,  as shown in Figure~\ref{fig:Fval4}. In particular,  if \(z^{f}\) lies in the \(\chi_{l}\)-character space then the lines from \(e_{1}\) and \(e_{2}\) are marked with \(\chi_{l}\).  Finally,  as \(L_{3,\beta}\) defeats the lines from \(e_{1}\) and \(e_{2}\) at \(v\), we have \(c > b\) by Remark~\ref{remark:knock-out}.  Hence \(\chi_{k}\neq \chi_{l}\).
 \qed

 \noindent\textsc{Case 3:} A vertex \(v\) of valency 5 or 6 (excluding three straight lines meeting at a point) defines a surface scroll blown-up in one or two points. 
  \begin{lemma}
 \label{lemma:val5or6}
 There are uniquely determined characters \(\chi_{k}\) and \(\chi_{l}\) which each mark a pair of lines meeting at \(v\).  The remaining line or pair of lines are marked with distinct characters. Mark the vertex \(v\) with \(\chi_{m} \! := \chi_{k}\otimes \chi_{l}\).  
 \end{lemma}
 \proof
 A vertex of valency 5 occurs only at the intersection point of a line \(L_{\alpha,\beta}\) from \(e_{\alpha}\) with a line \(L_{\gamma,\delta}\) from \(e_{\gamma}\).  We may assume that \(\alpha =3, \gamma = 1\), and that \(L_{3,\beta}\) is defeated so \(L_{1,\delta}\) extends.   This accounts for three lines meeting at \(v\);  the fourth and fifth lines are tesselating lines of a regular triangle \(T\) which is either a corner triangle from \(e_{1}\),  the meeting of champions triangle or a corner triangle from \(e_{2}\).  We illustrate the first case in Figure~\ref{fig:Fval5}: the lines \(L_{3,\beta}\) and \(L_{1,\delta}\) are cut out by \(x^{a}\! : \! y^{e}\) and \(y^{c}\! : \! z^{f}\) respectively,  while the tesselating lines extending from \(v\) into \(T\) are cut out by \(x^{d - r}\! : \! y^{b+r}z^{r}\),  for \(r\) satisfying the relations (\ref{eqn:prop3.1a}),  and by \(z^{h - i}\! : \! x^{i}y^{g+i}\),  for some \(g,h,i\) (\(i\neq 0\)).  

 The character \(\chi_{k}\) marking \(L_{1,\delta}\) marks two lines meeting at \(v\) as \(L_{1,\delta}\) passes through \(v\).  From (\ref{eqn:prop3.1a}) we have \(x^{d - r} = x^{a}\),  therefore a character \(\chi_{l}\) marks the lines cut out by both \(x^{a}\! : \! y^{e}\) and \(x^{d - r}\! : \! y^{b+r}z^{r}\),  so \(\chi_{l}\) also marks a pair of lines at \(v\).  Finally,  the character \(\chi_{j}\) marking the fifth line is neither \(\chi_{k}\) nor \(\chi_{l}\).  Indeed,  the relations (\ref{eqn:prop3.1a}) for \(T\) ensure that \(h-i > f\), so \(\chi_{j}\neq \chi_{k}\).  To show that \(\chi_{j} \neq \chi_{l}\),  write \(\tau\) for the basic triangle in \(T\) with two edges cut out by \(x^{d - r}\! : \! y^{b+r}z^{r}\) and \(z^{h - i}\! : \! x^{i}y^{g+i}\).  The corresponding toric variety is \(U_{\tau} = \Spec\; \C[\lambda,\mu,\nu]\),  where \(\lambda, \nu\) satisfy \(y^{b+r}z^{r} = \lambda x^{d-r}\) and \(x^{i}y^{g+i} = \nu z^{h-i}\) (see \cite[\S5]{Craw:htc}).  Thus,  both \(z^{h-i}\) and \(x^{d-r}\) lie in the basis of \(\owe_{Z}\) for the \(A\)-cluster \(Z\) defined by the origin \(\lambda=\mu=\nu=0\) in \(\C^{3}\cong U_{\tau}\).  It follows that \(z^{h-i}\) and \(x^{d-r}\) lie in different character spaces, so \(\chi_{j} \neq \chi_{l}\).  This proves the lemma when \(T\) is a corner triangle from \(e_{1}\).  

 When \(T\) is the meeting of champions or a corner triangle from \(e_{2}\),  the ratios cutting out the fourth and fifth lines are either \(x^{d - i}\! : \! y^{b+i}z^{i}\) with \(z^{h - k}\! : \! x^{g+k}y^{k}\),  or \(x^{d - i}\! : \! y^{i}z^{b+i}\) with \(z^{h - k}\! : \! x^{g+k}y^{k}\).  In each case the same argument applies so the lemma is established for a vertex of valency 5.  The case where the vertex has valency 6 is similar.
 \qed

 \noindent\textsc{Case 4:}  A vertex \(v\) at the intersection of three straight lines in \(\Sigma\) defines an exceptional del Pezzo surface of degree six,  denoted \(\dP\).  
 \begin{lemma}
 \label{lemma:val6}
 The monomials defining the pair of morphisms \(\dP\to \mathbb{P}^{2}\) lie in uniquely determined character spaces \(\chi_{l}\) and \(\chi_{m}\) satisfying 
 \begin{equation}
 \label{eqn:val6relation}
 \chi_{l} \otimes \chi_{m} = \chi_{i} \otimes \chi_{j} \otimes \chi_{k},
 \end{equation} 
 where \(\chi_{i}\), \(\chi_{j}\) and \(\chi_{k}\) mark the straight lines through the vertex \(v\) defining the del Pezzo surface.  Mark the vertex \(v\) with both \(\chi_{l}\) and \(\chi_{m}\).
 \end{lemma}
 \proof
 Three straight lines intersect at a vertex \(v\) in \(\Sigma\) only when three lines tesselating the same regular triangle intersect.  If \(v\) lies in a corner triangle then the three ratios listed in (\ref{eqn:prop3.2a}) which cut out the lines satisfy \(i+j+k=r\) (the case where the lines are cut out by the ratios (\ref{eqn:prop3.2b}) is similar).  The ratios determine a Segre embedding \(\mathbb{P}^{1}\times\mathbb{P}^{1}\times\mathbb{P}^{1}\rightarrow \mathbb{P}^{7}\) given by
 \begin{eqnarray*}
\lefteqn{(x^{d-i}y^{e-j}z^{f-k}\!:\!x^{d-i+k}y^{e-j+c+k}\!:\!x^{d-i+a+j}z^{f-k+j}\!:\!x^{d-i+a+j+k}y^{c+k}z^{j}\!:} \\
 & & y^{e-j+b+i}z^{f-k+i}\!:\! x^{k}y^{e-j+c+k+b+i}z^{i}\!:\!x^{a+j}y^{b+i}z^{f-k+i+j}\!:\!x^{a+j+k}y^{b+i+c+k}z^{i+j}).
 \end{eqnarray*}
 The del Pezzo \(\dP\subset \mathbb{P}^{6}\) is the intersection of the image of this map with the hyperplane \(x_{0} = x_{7}\),  where \(x_0,\dots ,x_7\) are coords on \(\mathbb{P}^{7}\).  Moreover, the maps \(\dP\rightarrow \mathbb{P}^{2}\) are the restriction of the projections \((x_0\!:\!x_2\!:\!x_3)\) and \((x_0\!:\!x_4\!:\!x_5)\) to \(\dP\).  After removing common factors and simplifying exponents using (\ref{eqn:prop3.1a}), these projections are
 \begin{equation}
 \label{eqn:mapstoP2}
 \left(y^{e-j}z^{i}\!:\!x^{a+j}z^{f-k}\!:\!x^{d-i}y^{c+k}\right)\quad\mbox{and}\quad\left(x^{d-i}z^{j}\!:\!y^{b+i}z^{f-k}\!:\!x^{k}y^{e-j}\right).
 \end{equation}
The required characters \(\chi_{l}\) and \(\chi_{m}\) are the common characters of the monomials defining these maps, i.e.,  the characters of say \(x^{a+j}z^{f-k}\) and \(x^{d-i}z^{j}\).  The product of this pair equals the product of \(x^{d-i}\),  \(z^{f-k}\) and \(x^{a+j}z^{j}\),  so (\ref{eqn:val6relation}) holds.
 \qed

 \begin{remark}
 \label{remark:passthrough}
 \emph{When two or more lines marked with the same character meet at a vertex \(v\) it is convenient to regard the lines as a single line \emph{passing through} \(v\).  Thus,  for example,  a valency 4 vertex \(v\) is the intersection point of two lines passing through \(v\).}
 \end{remark}

 \begin{example}
 \label{ex:11.1.2.8rr}
 \emph{Consider once again the cyclic quotient singularity of type \(\textstyle{\frac{1}{11}(1,2,8)}\).}  
 \begin{figure}[!ht]
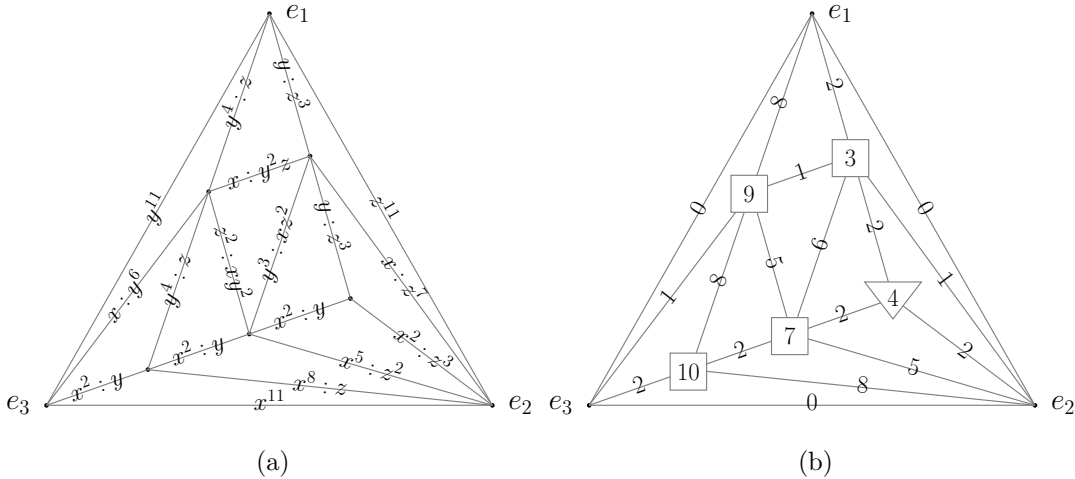

 \centering
 \mbox{\subfigure[]{\input{ratios.pspic}\label{fig:ratios}} \hspace{0.8cm} \subfigure[]{\input{11.1.2.8rr.pspic}\label{fig:11.1.2.8rr}}}
 \caption{(a) Ratios on lines in $\Sigma$; (b) Reid's recipe for $\frac{1}{11}(1,2,8)$}
 \end{figure}
 \emph{The ratios cutting out the lines in \(\Sigma\) are shown in Figure~\ref{fig:ratios}.  For \(\varepsilon\) a primitive 11th root of unity,  write \(\chi_{i} = \varepsilon^{i}\) (\(i = 0,\dots ,10\)) for the characters of \(A = \Z/11\).  The lines meeting at the vertex of valency 3 are marked with \(\chi_{2}\),  the common character space of the monomials \(x^{2},y,z^{3}\).  According to Lemma~\ref{lemma:val3} we mark the vertex of valency 3 with \(\chi_{4}\).  The characters \(\chi_{2}\) and \(\chi_{8}\) mark lines passing through (in the sense of Remark~\ref{remark:passthrough}) the vertex of valency 4 so,  by Lemma~\ref{lemma:val4},  we mark the vertex with \(\chi_{10}\).  The remaining vertices have valency 5,  so Lemma~\ref{lemma:val5or6} applies.  The result is shown in Figure~\ref{fig:11.1.2.8rr}.}
 \end{example}

 \begin{example}
 \label{ex:30.25.2.3rr}
 \emph{The fan \(\Sigma\) of \(\ahilb\) for the \(A\)-action \(\frac{1}{30}(25,2,3)\) is shown in Figure~\ref{fig:30.25.2.3rr}.  There are three regular triangles of side 2 to the left of the line from \(e_{1}\) to \(p\),  and two regular triangles of side 3 to the right.
 \begin{figure}[!ht]
 \begin{center}\input{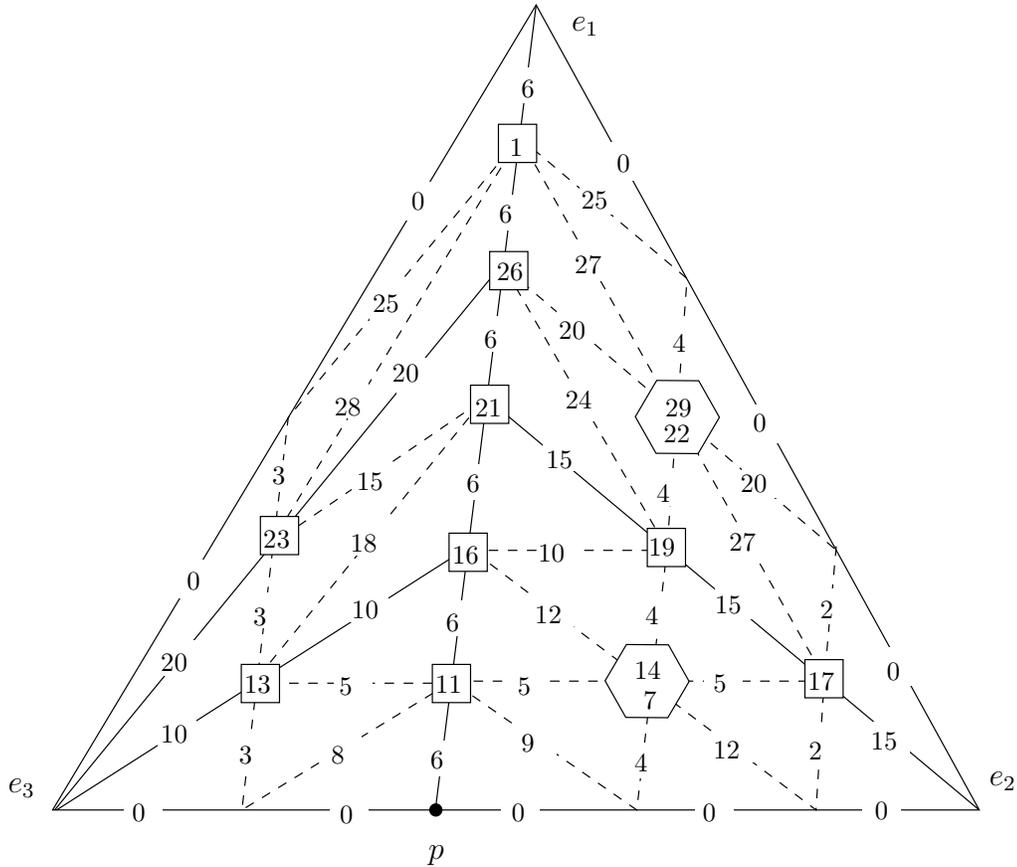}\end{center}
 \caption{Reid's recipe for $\frac{1}{30}(25,2,3)$}
 \label{fig:30.25.2.3rr}
 \end{figure}
 Every internal vertex has valency 5 or 6.  Most of the vertices are marked with a single character determined by Lemma~\ref{lemma:val5or6}.  However,  inside each regular triangle of side 3 is a vertex of valency 6 defining a del Pezzo surface \(\dP\), so each of these vertices is marked with a pair of characters.  For example,  \(\chi_{4}\), \(\chi_{5}\) and \(\chi_{12}\) mark the lines passing through one of the vertices.  The proof of Lemma~\ref{lemma:val6} reveals that the monomials defining the morphisms \(\dP\to\mathbb{P}^{2}\) lie in the \(\chi_{7}\) and \(\chi_{14}\) character spaces, so these characters mark the vertex.}
\end{example}

 \section{Every character appears once on \protect\(\Sigma\protect\)}
 \label{sec:meswc}
 It is not clear a priori from the construction of the previous section that different vertices are marked with different characters.  Nevertheless,  this is the case in Examples~\ref{ex:11.1.2.8rr} and~\ref{ex:30.25.2.3rr} where,  in addition,  every character of \(A\) marks either a line in \(\Sigma\) (possibly passing through several vertices in the sense of Remark~\ref{remark:passthrough}) or a unique vertex.   In this section we prove that this is the case for every finite Abelian subgroup \(A\subset \SL(3,\C)\).

\begin{figure}[!ht]
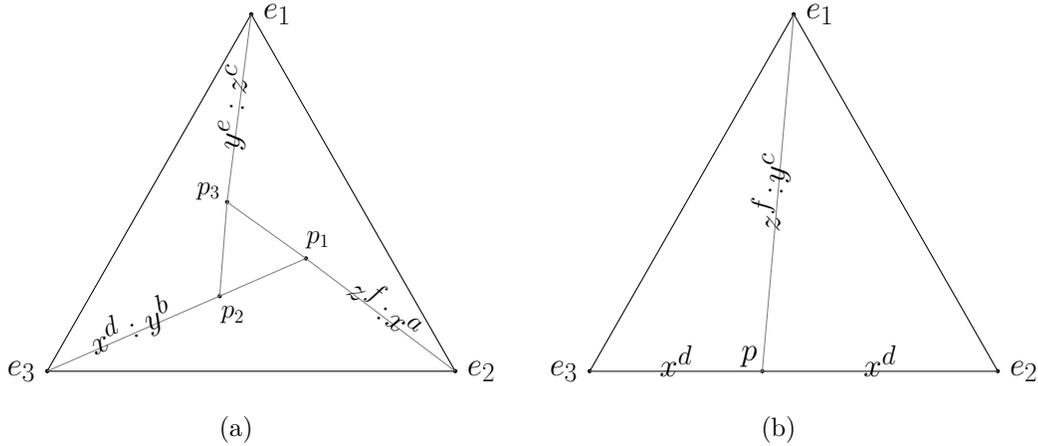

 \centering
 \mbox{\subfigure[]{\input{champ.pspic}\label{fig:champ}} \hspace{1.8cm} \subfigure[]{\input{long.pspic}\label{fig:long}}}
 \caption{Coarse subdivision: (a) meeting of champions; (b) long side}
  \label{fig:dichotomy}
 \end{figure}

 There is a dichotomy in the calculation of \(\Sigma\):  there is either a unique \emph{meeting of champions} or a unique \emph{long side} (see \cite[\S2.8.2]{Craw:htc}).   If a meeting of champions exists,  the champion lines subdivide \(\Sigma\) into four regions (see Figure~\ref{fig:champ}),  three if the champion has side zero or one if the meeting of champions is the whole of \(\Delta\).    
 Otherwise,  permuting \(x,y,z\) if necessary,  \(\Sigma\) is subdivided by a line from \(e_{1}\) which cuts the long side \(e_{2}e_{3}\) as in Figure~\ref{fig:long}  (there may be more than one line from \(e_1\) cutting the long side so this subdivision is not canonical).  In each case we produce a coarse subdivision of \(\Sigma\) into at most four regions which are themselves unions of regular triangles.  Each region,  apart from the interior triangle in Figure~\ref{fig:champ},  is a triangle with vertices \(e_{i}, p_{j}, e_{k}\);  in Figure~\ref{fig:long},  the point \(p = p_{j}\) lies on the edge \(e_{2}e_{3}\) cut out by \(x^d\).  Example~\ref{ex:11.1.2.8rr} contains a meeting of champions of side zero so is divided into three regions, Example~\ref{ex:30.25.2.3rr} has a long side and hence two regions.

 From now on we identify characters of \(A\) with monomials in the eigenspace of that character.  There is of course no canonical monomial for each character.  However,  we now show that the characters marking the points and lines in a regular corner triangle prefer a single monomial above all other choices (the meeting of champions triangle is different,  see Remark~\ref{rem:symmetry}).

 \begin{proposition}
 \label{prop:regions}
 The characters which mark the points and lines lying in the region \(e_{1}p_{2}e_{3}\) of Figure~\ref{fig:champ} can be represented by the monomials
 \begin{equation}
 x^{i}z^{j}\quad \mbox{for} \quad i = 0,\dots ,d\: ;  j = 0,\dots ,f. 
 \label{eqn:region}
 \end{equation}
 \end{proposition}

 \noindent By permuting \(x,y,z\) if necessary,  this proposition computes the characters marking any of the outer regions \(e_{i}p_{j}e_{k}\) in Figure~\ref{fig:champ} or \ref{fig:long}.  

\begin{lemma}
\label{lemma:chars_triangles}
The characters which mark the regular triangles of side \(r\) in Figure~\ref{fig:reg_triangles} can be represented by the monomials
\begin{align}
z^{f-k} \quad & \mbox{and}\quad x^{d-i}z^{f-k}\quad \mbox{for}\quad i,k = 0,\dots ,r \quad & \mbox{in Case a,} \label{eqn:charsa} \\
x^{i}z^{f-k} \quad & \mbox{and}\quad x^{d-i}z^{c+k}\quad  \mbox{for} \quad 0\leq i+k \leq r & \mbox{in Case b.}\label{eqn:charsb}
\end{align}
\end{lemma}
\noindent\emph{Proof of Proposition~\ref{prop:regions},  assuming the lemma.\ }
Starting from the edge \(e_{1}e_{3}\),  run an MMP (see \cite[\S2.7]{Craw:htc}) which eats all regular triangles inside the region \(e_{1}p_{2}e_{3}\) of Figure~\ref{fig:champ}.  We prove the proposition by induction on the number of contractions in the MMP.  If the MMP consists of a single contraction then the region itself is a regular corner triangle from \(e_{3}\),  shown in Figure~\ref{fig:corner_triangle}.  The ratio \(x^{a}\!:\!y^{e}\) cutting out the edge \(e_{1}e_{3}\) is simply \(y^{e}\),  hence \(a = 0\).  Since \(d - a = f = r\) holds by (\ref{eqn:prop3.1a}),  substitute \(d = f = r\) into the list (\ref{eqn:charsa}) of characters marking a corner triangle to see that the proposition holds in this case.

 Suppose now that we have performed an MMP that has eaten all regular triangles in a region with vertices \(e_{1}qe_{3}\) where the lines \(e_{1}q\) and \(e_{3}q\) are cut out by the ratios \(z^{f}\!:\!y^{c}\) and \(x^a\!:\!y^e\) respectively (see Figure~\ref{fig:corner_triangle}).  We assume by induction that the characters marking the union of regular triangles inside this region are 
\begin{equation}
\label{eqn:charsinduction}
x^{i}z^{j}\quad \mbox{for} \quad i = 0,\dots ,a;  j = 0,\dots ,f.
\end{equation}
If the next contraction of the MMP eats a corner triangle from \(e_{3}\),  then the line \(e_{1}q\) extends to a lattice point \(p\),  and the line \(e_{3}p\) has ratio \(x^d\!:\!y^b\) say,  as shown in Figure~\ref{fig:corner_triangle}.  The characters which mark the new corner triangle are listed in (\ref{eqn:charsa}).  The region \(e_{1}p_{2}e_{3}\) of Figure~\ref{fig:champ} is therefore marked with the union of characters (\ref{eqn:charsa}) and (\ref{eqn:charsinduction});  namely \(x^{i}z^{j}\) for \(i = 0,\dots d;\;  j = 0,\dots f\) as required.   The case where the final triangle is from \(e_{1}\) is similar. \qed

\noindent\emph{Proof of Lemma~\ref{lemma:chars_triangles}.\ } For Case a,  the triangle is eaten by an MMP from the side \(e_{1}e_{3}\) so we choose to represent the characters marking this triangle by monomials in \(x,z\).  From (\ref{eqn:prop3.2a}),  the characters which mark the tesselating lines of the triangle are 
\begin{equation}
\label{eqn:linechars}
x^{d-i}; \quad x^{a+j}z^{j}; \quad z^{f-k}\quad \mbox{for} \quad i,j,k = 0,\dots r-1.
\end{equation}
The vertices along the edges of the triangle are marked with the characters 
\begin{equation}
\label{eqn:pointchars}
x^{a}z^{f-k}; \quad x^{d}z^{f-k};\quad x^{d-i}z^{f}\quad\mbox{for}\quad i,k = 0,\dots r-1.
\end{equation}
Indeed,  the edges emanating from \(e_{3}\) are marked with \(x^{a}\) and \(x^{d}\).  A tesselating line marked with \(z^{f-k}\) (for \(k = 0,\dots r-1\)) passes through every vertex on both of these edges and hence,  by Lemmas~\ref{lemma:val4} and~\ref{lemma:val5or6},  these vertices are marked with \(x^{a}z^{f-k}\) and \(x^{d}z^{f-k}\).  Similarly,  the tesselating lines \(x^{d-i}\) cross the edge from \(e_{1}\) marked with \(z^{f}\),  so \(x^{d-i}z^{f}\) mark the vertices along this edge.  Finally,  from the proof of Lemma~\ref{lemma:val6},  we know that the characters 
 \begin{equation}
 \label{eqn:valency6_a}
x^{a+j}z^{f-k} \mbox{ and } x^{d-i}z^{j}\quad \mbox{for } i,j,k = 1,\dots r-1 \; \mbox{such that}\; i+j+k=r
 \end{equation}
mark the internal vertices in the tesselation of the regular triangle of Figure~\ref{fig:corner_triangle}.  As a result,  the union of the characters listed in (\ref{eqn:linechars}), (\ref{eqn:pointchars}) and (\ref{eqn:valency6_a}) mark the points and lines of the regular triangle of Figure~\ref{fig:corner_triangle}.  It is an easy combinatorial exercise to see that the union of these characters is equal to the list (\ref{eqn:charsa}) as claimed.  This proves Case a of Lemma~\ref{lemma:chars_triangles}.  

 \indent To prove Case b,  one proves similarly that the characters \(x^{d-i}\), \(x^{j}z^{c+j}\) and \(z^{f-k}\) for \(i,j,k = 0,\dots ,r-1\) mark the lines of the regular tesselation;  the characters \(x^{d}z^{f-k}\), \(x^{a+j}z^{c+j}\) and \(x^{d-i}z^{c}\) for \(i,j,k = 0,\dots ,r-1\) mark the vertices along the edges of the triangle;  and the characters \(x^{j}z^{f-k}\) and \(x^{d-i}z^{c+j}\) for \(i,j,k = 1,\dots r-1\) such that \(i+j+k=r\) mark the internal vertices of the triangle.  As with Case a,  the union of these characters is equal to the list (\ref{eqn:charsb}) as claimed.\qed

 \begin{remark} 
 \label{rem:symmetry}
 \emph{There is symmetry in Lemma~\ref{lemma:chars_triangles} Case b.  The characters were listed in terms of \(x,z\),  but equally can be written in \(x,y\) or \(y,z\) using the relations (\ref{eqn:prop3.2b}).  This doesn't alter the character,  because the ratios in (\ref{eqn:prop3.2b}) are \(A\)-invariant.  In short, the characters marking strata in the interior triangle in Figure~\ref{fig:champ} do not prefer a single monomial over all others.}
 \end{remark}

\begin{figure}[!ht]
 \begin{center}\input{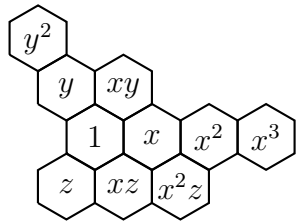}\end{center}
 \caption{The McKay quiver as a tesselation by regular hexagons}
 \label{fig:quiver_xyz}
 \end{figure}

 The condition \(A\subset \SL(3,\C)\) ensures that \(xyz\) is \(A\)-invariant, so a monomial lies in the same character space as its image in \(\C[x,y,z]/xyz\).  Following Reid~\cite[\S7]{Reid:mc},  the monomials in \(\C[x,y,z]/xyz\) are represented as a tesselation of the plane by regular hexagons,  part of which is shown in Figure~\ref{fig:quiver_xyz}. 
 
 This is the universal cover of the \emph{McKay quiver},  where the arrows in the three principal directions are `multiply by \(x,y\) or \(z\)'.  Some power of each monomial \(x,y\) and \(z\) is \(A\)-invariant so the tesselation is periodic,  and we say that any connected region in the quiver in one-to-one correspondence with the characters of \(A\) is a \emph{fundamental domain}.

\begin{figure}[!ht]
 \begin{center}\input{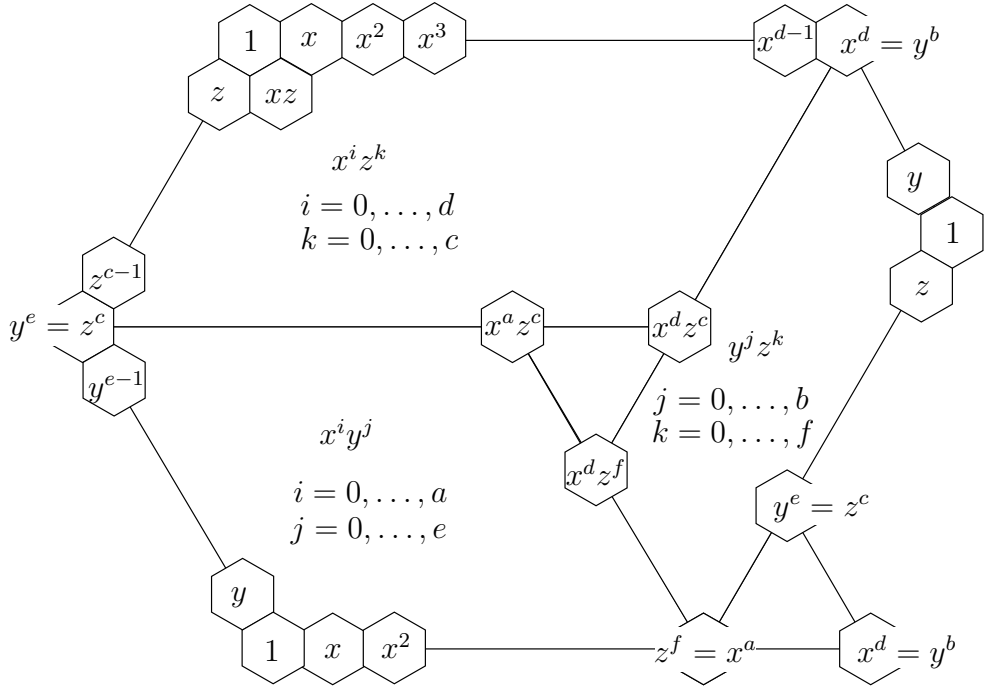}\end{center}
 \caption{Three parallelograms and two triangles in the McKay quiver}
 \label{fig:pf_main_thm}
 \end{figure}

 \begin{proposition} 
 \label{prop:one-to-one}
 The characters marking the points and lines in \(\Sigma\) form a fundamental domain in the McKay quiver (assuming that a character which marks a line passing through a vertex in the sense of Remark~\ref{remark:passthrough} is recorded only once on the quiver). 
 \end{proposition}
 \proof The coarse subdivision of \(\Sigma\) is one of the two types shown in Figure~\ref{fig:dichotomy}.  Beginning with Case a,  plot the characters which mark each region on the McKay quiver.  The characters marking the three outer regions in the subdivision form parallelograms,  by Proposition~\ref{prop:regions},  and the characters marking the meeting of champions form a pair of triangles,  by Lemma~\ref{lemma:chars_triangles} Case b.  The parallelograms and triangles intersect along characters \(x^{i}y^{e}=x^{i}z^{c}\) (\(0\leq i \leq d\)),  \(y^{j}z^{f} = y^{j}x^{a}\) (\(0\leq j \leq e\)) and  \(x^{d}z^{k}=y^{b}z^{k}\) (\(0\leq k \leq f\)) marking the vertices on the champion lines, as shown in Figure~\ref{fig:pf_main_thm}.
The union of these regions is slightly larger than a fundamental domain in the quiver.  However,  the characters \(x^{i}, y^{j}, z^{k}\) around the edge of the shape in Figure~\ref{fig:pf_main_thm} mark tesselating lines in different regions and have been plotted more than once.  In each case,  the lines marked with the same monomial pass through a vertex of valency 4,5 or 6 (c.f. Remark~\ref{remark:passthrough}),  thereby passing from one region to another.  The assumption in the statement of the theorem enables us to identify in pairs (the trivial character 1 appears three times) the monomials around the outside of the shape of Figure~\ref{fig:pf_main_thm},  leaving a fundamental domain.

 Otherwise,  the subdivision is from Figure~\ref{fig:long}.  The characters marking the two regions are \(x^{i}z^{k}\) and \(x^{i}y^{j}\) for \(i = 0,\dots d\),  \(j = 0,\dots c\) and \(k = 0,\dots ,f\). These parallelograms intersect along the characters \(x^{i}y^{c} = x^{i}z^{f}\) (\(0\leq i \leq d\)) marking the vertices on the line of intersection of the regions.  Since \(x^{d}\) is \(A\)-invariant,  identify the characters \(y^{j}\) and \(x^{d}y^{j}\) pairwise,   and similarly \(z^{k}\) and \(x^{d}z^{k}\).  Identify in pairs the two collections \(1,x,\dots,x^{d}\) marking the lines passing from one region to another,  leaving a fundamental domain.\qed

 \begin{remark}
 \emph{The coarse subdivision of Figure~\ref{fig:long} is not canonical,  different subdivisions vary by a corner triangle \(T\) of side \(r\) from \(e_{1}\) whose sides extend from \(e_{1}\) to the long side.  If the long side is cut out by the monomial \(x^{f}\) and the other sides of \(T\) are cut out by \(y^{a}\!:\!z^{e}\) and \(y^{d}\!:\!z^{b}\) say,  where the relations (\ref{eqn:prop3.2a}) hold,  then the characters which mark \(T\) are \(y^{d-i}x^{f-k}\) for \(i,k = 0,\dots ,r\) by Lemma~\ref{lemma:chars_triangles}.  This is a parallelogram with sides of length \(r\), so choosing a different coarse subdivision causes the pair of parallelograms in the above proof to be translated.  In particular,  the overall result is unchanged.}
 \end{remark}

  \begin{corollary} 
 \label{coro:chars}
 Every nontrivial character of \(A\) appears once on \(\Sigma\) as either
 \begin{enumerate}
 \item[\one] a character \(\chi_{i}\) marking a line (possibly passing through several vertices in the sense of Remark~\ref{remark:passthrough}); or
 \item[\two] a character \(\chi_{m}\) marking a vertex; or
 \item[\three] the second character \(\chi_{l}\) marking the intersection of three straight lines.
 \end{enumerate}
 \end{corollary}

 \section{Tautological line bundles on \ahilb}
 \label{sec:tlb} 
 For a finite Abelian subgroup \(A\subset \SL(3,\C)\),  write  \(\pi \colon \C^{3}\rightarrow X = \C^{3}/A\) for the quotient, \(Y = \ahilb\) for Nakamura's \(A\)-Hilbert scheme and \(\varphi\colon Y \to X\) for the crepant resolution.  Since \(Y\) is a fine moduli space of subschemes \(Z\subset \C^{3}\) there is a universal subscheme \(\mathcal{Z}\subset Y\times \C^{3}\) fitting into the diagram

 \begin{equation*}\label{maindiag}
 \setlength{\unitlength}{36pt}
 \begin{picture}(2.2,2)(0,0)
 \put(1,0){\object{$X$}}
 \put(0,1){\object{$Y$}}
 \put(2,1){\object{$\C^{3}$}}
 \put(1,2){\object{$\mathcal{Z}$}}
 \put(0.25,0.75){\vector(1,-1){0.5}}
 \put(0.5,0.5){\nwlabel{$\varphi$}}
 \put(1.75,0.75){\vector(-1,-1){0.5}}
 \put(1.5,0.5){\swlabel{$\pi$}}
 \put(0.75,1.75){\vector(-1,-1){0.5}}
 \put(0.5,1.5){\nelabel{$p$}}
 \put(1.25,1.75){\vector(1,-1){0.5}}
 \put(1.5,1.5){\selabel{$q$}}
 \end{picture}
\end{equation*}
in which \(\pi\) and \(p\) are finite,  \(\varphi\) and \(q\) are birational and \(p\) is flat.  The sheaf \(\mathcal{R} \! := p_{*}\owe_{\mathcal{Z}}\) is locally free since \(p\) is finite and flat.  Write \(Z(y)\subset \C^{3}\) for the \(A\)-cluster corresponding to a point \(y\in Y\).  By the definition of \(A\)-cluster (see \S\ref{sec:intro}),  the \(A\)-module \(H^{0}(Z(y),\owe_{Z(y)})\) is the regular representation of \(A\),  therefore \(\dimn H^{0}(Z(y),\owe_{Z(y)}) \) is constant,  namely the order of the group \(A\).  A theorem of Grothendieck on (higher) direct images of coherent sheaves under proper morphisms (see Mumford~\cite[II.5,Corollary 2]{Mumford:av}) establishes that the fibre of \(\mathcal{R}\) over \(y\) is \(H^{0}(Z(y),\owe_{Z(y)})\).  In particular,  the rank of \(\mathcal{R}\) is equal to the order of the group \(A\). The decomposition of the regular representation into irreducible submodules induces the decomposition
 \[
 \mathcal{R} = \bigoplus_{k}\mathcal{R}_{k}\otimes \rho_{k}\quad \mbox{for }\mathcal{R}_{k} = \Hom_{A}(\rho_{k},\mathcal{R}),
 \]
 where the sum runs over all irreducible representations \(\rho_{k}\) of \(A\).  The sheaves \(\mathcal{R}_{k}\) are direct summands of a locally free sheaf so are themselves locally free of rank \(\dim \rho_{k} = 1\)  We call \(\mathcal{R}_{k}\) the \emph{tautological line bundle} on \(Y\) associated to the irreducible representation \(\rho_{k}\) of \(A\).

 Our calculation of the coordinates on the subvarieties \(U_{\tau} \cong \C^{3}\) covering \(Y\) enables us to compute an explicit basis of the fibres of \(\mathcal{R}\) over \(U_{\tau}\).  First,  write down the coordinates \(\xi,\eta,\zeta\) (or \(\lambda,\mu,\nu\)) on \(U_{\tau}\) given in one of the lists (\ref{eqn:coordsupa}) to (\ref{eqn:coordsdownb}).  In each case,  the origin \(0\in \C^{3}\cong U_{\tau}\) defines an \(A\)-cluster \(Z_{\tau}\) with defining ideal \(I_{\tau}\) generated by monomials.  It is easy to write down a \(\C\)-basis for \(H^{0}(Z_{\tau},\owe_{Z_{\tau}}) = \C[x,y,z]/I_{\tau}\),  namely the set \(\Gamma_{\tau}\) of monomials in \(\C[x,y,z]\smallsetminus I_{\tau}\).  This set is called an \emph{\(A\)-graph} by Nakamura.  Every monomial \(m\in \C[x,y,z]\) lies in a well defined character space of the \(A\)-action,  giving rise to a map \(\wt\colon \Gamma_{\tau} \to A^{\vee}\) from an \(A\)-graph to the character group of \(A\). This map is one-to-one because \(H^{0}(Z_{\tau},\owe_{Z_{\tau}})\) is the regular representation.  Nakamura~\cite[Lemma~2.3(ii)]{Nakamura:ago} proves that the set \(\Gamma_{\tau}\) forms a basis of the coordinate ring \(H^{0}(Z(y),\owe_{Z(y)})\) for every point \(y\) in the affine chart \(U_{\tau}\) (see also \cite[\S4]{Craw:htc}).    Thus we have shown the following well known fact:

 \begin{proposition}
 The monomials in the \(A\)-graph \(\Gamma_{\tau}\) form a basis of the fibres of \(\mathcal{R}\) over every point of \(U_{\tau}\).
 \label{prop:Agraph}
 \end{proposition}

\begin{figure}[!ht]
 \centering
 \input{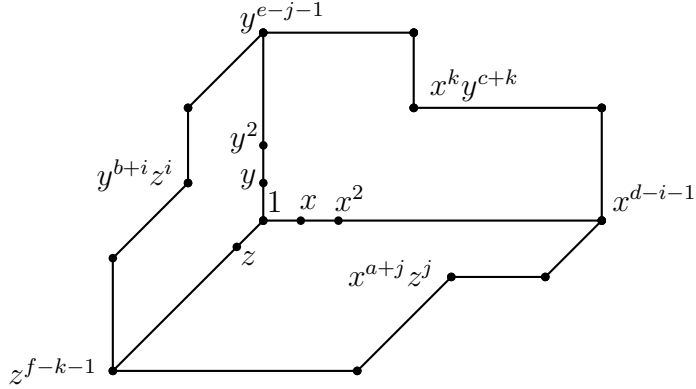}
 \caption{The $A$-graph $\Gamma_{\tau}$ whose monomials base the fibres of $\mathcal{R}$ over $U_{\tau}$.}
 \label{fig:tripod}
 \end{figure}

 To illustrate this,  consider the chart \(U_{\tau}\subset Y\) with coordinates \(\xi,\eta,\zeta\) as in (\ref{eqn:coordsupa}).  The point \(\xi = \eta = \zeta = 0\) corresponds to the \(A\)-cluster with ideal 
 \[
 I_{\tau} = \langle x^{d-i},y^{e-j},z^{f-k},y^{b+i+1}z^{i+1},x^{a+j+1}z^{j+1},x^{k+1}y^{c+k+1},xyz\rangle.
 \] 
 The \(A\)-graph \(\Gamma_{\tau}\) consisting of the monomials lying in \(\C[x,y,z]\smallsetminus I_{\tau}\) is shown in Figure~\ref{fig:tripod}.  
 As before,  monomials in \(\C[x,y,z]/xyz\) are drawn on a planar graph,  but to save space we represent monomials by dots rather than hexagons.  For each \(\rho_{k}\), the generator \(r_{k,\tau}\) of the tautological line bundle \(\mathcal{R}_{k}\) over $U_{\tau}$ can be read off directly from \(\Gamma_{\tau}\): if \(\chi_{k}\in A^{\vee}\) denotes the character of the representation \(\rho_{k}\) then \(r_{k,\tau}\) is the unique monomial \(m\in \Gamma_{\tau}\) with \(\wt(m) = \chi_{k}\).

 \begin{remark}
 \emph{For finite \(G\subset \SL(2,\C)\),  Gonzalez-Sprinberg and Verdier~\cite{Gonzalez-Sprinberg:cgm} calculated the generators of each \(\mathcal{R}_{k}\) over an open cover of \(Y\) by computing a resolution of \(\mathcal{R}_{k}\).  The new approach via Nakamura's \(A\)-graphs \(\Gamma_{\tau}\) implemented here provides the generators of \(\mathcal{R}_{k}\) on the open cover of \(Y\) without the need for explicit resolutions. 
}
 \end{remark}

 The next lemma lists elementary facts about the generators \(r_{k,\tau}\).  Given \(\rho_{k}\),  let \(\mathfrak{m}\in \C[x,y,z]\) be such that \(r_{k,T} = \mathfrak{m}\) for some basic triangle \(T\). Set
 \[
 \Region(k,\mathfrak{m})\!:= \{p\in \Delta \st \exists \mbox{ basic triangle }\tau \mbox{ containing } p \mbox{ such that } r_{k,\tau} = \mathfrak{m}\}.
 \]
 This subset of the junior simplex \(\Delta\) is by construction a union of triangles. 

 \begin{lemma}
 \label{lemma:facts}
 \begin{enumerate}
 \item[\one] \(\Region(k,\mathfrak{m})\) is a convex subset of \(\Delta\).
 \item[\two] \(e_{1}\notin\Region(k,\mathfrak{m}) \iff \mathfrak{m}\) is divisible by \(x\).  The same statement holds if \(e_{2}\) (resp.\ \(e_{3}\)) replaces \(e_{1}\)  and \(y\) (resp.\ \(z\)) replaces \(x\).
 \item[\three] Let \(v\) be a vertex of triangles \(T, T'\).  If \(r_{k,T} = x^{\alpha}z^{\gamma}\) and \(r_{k,T'} = y^{\beta}z^{\delta}\) then \(z^{\min\{\gamma,\delta\}}\) divides \(r_{k,\tau}\) for all triangles \(\tau\) whose interior intersects \(e_{1}ve_{2}\).
 \end{enumerate}
 \end{lemma}
 \proof 
 Write \(\mathcal{R}_{k} \cong \owe_{Y}(D)\) for some divisor \(D\) on \(Y\).  Toric geometry defines \(D\) by specifying an element \(m_{\tau}\in M\) for each \(\tau\in \Sigma\),  defining divisors \(\divisor(m_{\tau}^{-1})\) on \(U_{\tau}\subset Y\) so that \(m_{\tau}\) generates \(\owe_{Y}(D)\vert_{U_{\tau}}\).  Now,  \(r_{k,\tau}\) generates \(\mathcal{R}_{k}\vert_{U_{\tau}}\) so,  accounting for linear equivalence,  there exists a Laurent monomial \(f\) such that \(r_{k,\tau} = f\cdot m_{\tau}\) (if \(f = 1\) then \(r_{k,\tau}\) is \(A\)-invariant,  but this is false unless \(\rho_{k}\) is trivial).  The line bundle \(\mathcal{R}_{k}\) is generated by its global sections so the piecewise linear function \(\psi_{D}\colon \vert \Sigma\vert \rightarrow \mathbb{R}\) defined by \(\psi_{D}(v) = \langle m_{\tau},v\rangle\) for \(v \in \tau\) is convex.  Since \(f\) is fixed,  it follows that \(\psi_{k}\colon \vert \Sigma\vert \rightarrow \mathbb{R}\) defined by \(\psi_{k}(v) = \langle r_{k,\tau},v\rangle\) for \(v \in \tau\) is also convex.  This proves part \one.  

 For any triangle \(T\) with vertex \(e_{1}\),  one of the coordinates on the affine variety \(U_{T}\) is \(x/y^{j}z^{k}\) for some \(j,k\in \Z_{\geq 0}\), so \(x\) cannot divide any monomial \(r_{k,T}\) in \(\Gamma_{T}\).  In particular,  when \(\mathfrak{m}\) is divisible by \(x\) we have \(e_{1}\notin\Region(k,\mathfrak{m})\).  Conversely,  if \(\mathfrak{m} = r_{k,\tau}\) is not divisible by \(x\) then \(\mathfrak{m} = y^{\beta}z^{\gamma}\) for \(\beta,\gamma\in \Z_{\geq 0}\).  Then \(\psi_{k}(v) = \langle y^{\beta}z^{\gamma},v\rangle\) for \(v\in \tau\) and \(\langle y^{\beta}z^{\gamma},e_{1}\rangle = 0\).  But \(\psi_{k}(e_{1}) = 0\) because \(x\) cannot divide \(r_{k,T}\) for any \(T\) with vertex \(e_{1}\),  so in fact \(\psi_{k}(v) = \langle y^{\beta}z^{\gamma},v\rangle\) for \(v = e_{1}\) and for all \(v\in \tau\).  Convexity of \(\psi_{k}\) ensures \(\psi_{k}(v) \leq \langle y^{\beta}z^{\gamma},v\rangle\) for all \(v\in \Sigma\),  so piecewise linearity gives \(\psi_{k}(v) = \langle y^{\beta}z^{\gamma},v\rangle\) for all \(v\) on straight lines in \(\Delta\) joining \(e_{1}\) to points of \(\tau\).  In particular,  \(r_{k,T} = y^{\beta}z^{\gamma}\) for some \(T\) with vertex \(e_{1}\),  so \(e_{1}\in\Region(k,\mathfrak{m})\).  This completes \two.  

 Finally,  suppose there exists \(\tau\) such that \(z^{\min\{\gamma,\delta\}}\) doesn't divide \(r_{k,\tau}\) and an interior point \(w\) of \(\tau\) lies inside \(e_{1}ve_{2}\).  By exchanging \(x\) and \(y\) if necessary,  we may assume \(r_{k,\tau} = x^{a}z^{c}\) with \(c<\min\{\gamma,\delta\}\).  If \(a \leq \alpha\) then any \(A\)-graph containing \(x^{\alpha}z^{\gamma}\) also contains \(x^{a}z^{c}\),  but \(\wt(x^{\alpha}z^{\gamma}) = \chi_{k} = \wt(x^{a}z^{c})\),  so in fact \(a > \alpha\).  Now,  consider \(z^{\gamma - c}/x^{a-\alpha} = x^{\alpha}z^{\gamma}/x^{a}z^{c}\in M = N^{\vee}\).  The plane \((z^{\gamma - c}/x^{a-\alpha})^{\perp}\subset N\otimes \mathbb{R}\) defines a line \(l = (z^{\gamma - c}/x^{a-\alpha})^{\perp}\cap \Delta\) in the junior simplex \(\Delta\) passing through \(e_{2}\).  Since \(r_{k,T} = x^{\alpha}z^{\gamma}\) and \(v\in T\) we have \(\langle z^{\gamma - c}/x^{a-\alpha}, v\rangle  = \langle x^{\alpha}z^{\gamma}, v\rangle - \langle x^{a}z^{c}, v\rangle \leq 0\).  Similarly,  \(r_{k,\tau} = x^{a}z^{c}\) and \(w\) lies strictly inside \(\tau\) so \(\langle z^{\gamma - c}/x^{a-\alpha}, w\rangle = \langle x^{\alpha}z^{\gamma}, w\rangle - \langle x^{a}z^{c}, w\rangle > 0\).  Clearly \(\langle z^{\gamma - c}/x^{a-\alpha},e_{3}\rangle > 0\).  As a result, \(e_{3}\) and \(w\) lie on the same side of the line \(l\) through \(e_{2}\) while \(v\) lies either on the opposite side of \(l\) or on the line itself.  Either way, \(w\) cannot lie in the region \(e_{1}ve_{2}\),  a contradiction. \qed

 \begin{example}
 \label{ex:newtpoly}
 \emph{The fan \(\Sigma\) of \(\ahilb\) for the singularity \(\textstyle{\frac{1}{11}(1,2,8)}\) is shown in Figure~\ref{fig:11.1.2.8rr}.  The generators \(r_{3,\tau}\) of \(\mathcal{R}_{3}\vert_{U_{\tau}}\) are drawn on Figure~\ref{fig:R3generators};  note that \(r_{3,\tau} = xy\) for four triangles \(\tau\in\Sigma\).  Six convex regions \(\Region(3,\mathfrak{m})\) partition \(\Delta\) for certain \(\mathfrak{m}\) with \(\wt(\mathfrak{m}) = \chi_{3}\).  Parts \two\ and \three\ of  Lemma~\ref{lemma:facts} are essentially obvious when you consider the relative positions of the \(\mathfrak{m}\) in the Newton polygon shown in Figure~\ref{fig:R3newtpoly}.  Observe that \(y^{2}z^{4}\) and \(yz^{7}\) don't generate \(\mathcal{R}_{3}\) on any open set:  the ratios of consecutive monomials \(y^{3}z/y^{2}z^{4}\),  \(y^{2}z^{4}/yz^{7}\) and \(yz^{7}/z^{10}\) coincide so \(\mathcal{R}_{3}\) has degree three on one of the curves parametrised by \(y\! : \! z^{3}\) (see Lemma~\ref{lemma:deg1} for more on this point). 
 \begin{figure}[!ht]
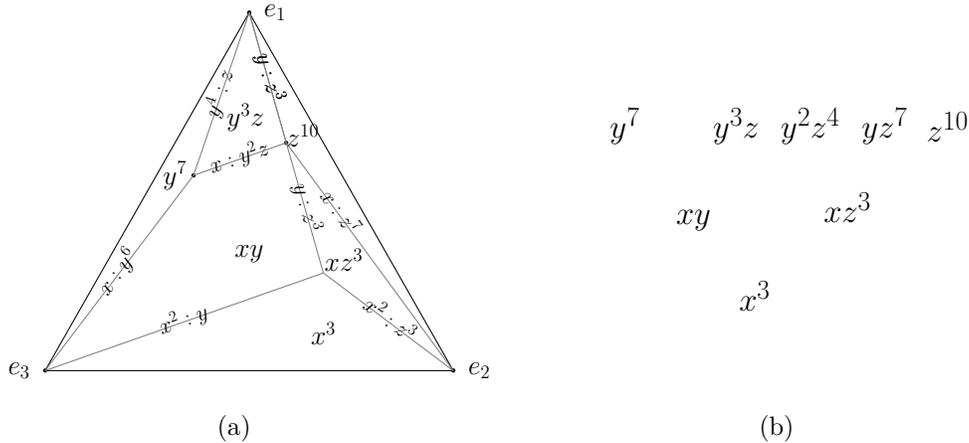

 \centering
 \mbox{\subfigure[]{\input{R3generators.pspic}\label{fig:R3generators}} \hspace{1.8cm} \subfigure[]{\input{R3newtpoly.pspic}\label{fig:R3newtpoly}}}
 \caption{(a) Generators $r_{3,\tau}$ of $\mathcal{R}_{3}$; (b) the Newton polygon}
 \end{figure}}
 \end{example}

 \section{Relations between tautological line bundles}
 \label{sec:rip}

 For every compact exceptional surface of the crepant resolution \(\varphi\colon Y\to X\),  there is a relation in \(\Pic(Y)\) between tautological line bundles \(\mathcal{R}_{k}\) on \(Y\).  These relations can be stated using Reid's recipe (see \S3):


 \begin{theorem} 
 \label{thm:bundles}
 The following relations hold in \(\Pic(Y)\):
 \begin{enumerate}
 \item \(\mathcal{R}_{m} = \mathcal{R}_{k}\otimes\mathcal{R}_{k}\) when \(\chi_{m} = \chi_{k}\otimes \chi_{k}\) marks a vertex  \(v\) of valency 3.
 \item \(\mathcal{R}_{m} = \mathcal{R}_{k}\otimes\mathcal{R}_{l}\) when \(\chi_{m} = \chi_{k}\otimes \chi_{l}\) marks a vertex \(v\) of valency 4. 
 \item \(\mathcal{R}_{m} = \mathcal{R}_{k}\otimes\mathcal{R}_{l}\) when \(\chi_{m} = \chi_{k}\otimes \chi_{l}\) marks a vertex \(v\) of valency 5 or 6.
 \item \(\mathcal{R}_{l}\otimes\mathcal{R}_{m} = \mathcal{R}_{i}\otimes\mathcal{R}_{j}\otimes\mathcal{R}_{k}\) when the pair of characters \(\chi_{l}\) and \(\chi_{m}\) satisfying  \(\chi_{l} \otimes \chi_{m} = \chi_{i} \otimes \chi_{j} \otimes \chi_{k}\) mark the intersection point \(v\) of three straight lines.
 \end{enumerate}
 \end{theorem}

 \proof To establish a relation of the form \(\mathcal{R}_{m} = \mathcal{R}_{k}\otimes\mathcal{R}_{l}\) we prove that  \(r_{m,\tau} = r_{k,\tau}\cdot r_{l\tau}\) on every triangle \(\tau\in \Sigma\). We proceed case by case as in \S\ref{sec:df}:

 \medskip

 \noindent\textsc{Case 1:} Write \(v\) for the vertex of valency 3 marked with \(\chi_{m} = \chi_{k}\otimes \chi_{k}\) from Lemma~\ref{lemma:val3}.  In the notation of Lemma~\ref{lemma:val3},  write \(T\) for the basic triangle in \(\Sigma\) with \(v\) as a vertex and ratios \(x^{d}\! : \! y^{b}\),  \(y^{b}\! : \! z^{f}\) cutting out two edges.  By permuting \(x,y,z\) if necessary,  assume that \(T\) lies in a regular corner triangle of side \(r\) from \(e_{3}\) as shown in Figure~\ref{fig:corner_triangle} (with \(b=c\)).  The third edge of \(T\) is cut out by \(y^{e-(r-1)}\! : \! x^{d+(r-1)}z^{r-1}\) so the coordinates on \(U_{T}\) are \(\xi = x^{d}/y^{b}\),  \(\eta = y^{e-(r-1)}/x^{d+(r-1)}z^{r-1}\) and \(\zeta = z^{f}/y^{b}\).  Calculating \(\Gamma_{T}\) using the method introduced in \S\ref{sec:tlb} shows that both \(y^{b}\) and \(y^{e-r}\) lie in \(\Gamma_{T}\)  Since \(\wt(y^{b}) = \chi_{k}\) we have \(r_{k,T} = y^{b}\).  Also,  \(e - r = 2b\) by (\ref{eqn:prop3.1a}) so \(\wt(y^{e-r}) = \wt(y^{2b}) = \chi_{k}\otimes \chi_{k} = \chi_{m}\),  hence \(r_{m,T} = y^{2b}\). Note in passing that \(y^{2b}\) lies in the socle of \(\Gamma_{T}\).

 We now claim that for every triangle \(\tau\) in the region \(e_{1}ve_{3}\) of \(\Sigma\) we have \(r_{k,\tau} = y^{b}\) and \(r_{m,\tau} = y^{2b}\) so that \(r_{m,\tau} = r_{k,\tau}\cdot r_{k,\tau}\).  Indeed,  \(\mathfrak{m} = y^{b}\) is divisible by neither \(x\) nor \(z\) so Lemma~\ref{lemma:facts}\two\ shows that both \(e_{1}\) and \(e_{3}\) lie in \(\Region(k,y^{b})\).  The vertex \(v\in T\) also lies in this set since \(r_{k,T} = y^{b}\) by the above.  It follows that the entire region \(e_{1}ve_{3}\) lies in \(\Region(k,y^{b})\) by convexity,  thereby proving the claim for \(r_{k,\tau} = y^{b}\).  The proof for \(r_{m,\tau} = y^{2b}\) is identical.  Symmetrically,  we see that \(r_{k,\tau} = x^{a}\) and \(r_{m,\tau} = x^{2a}\) for \(\tau\) in \(e_{2}ve_{3}\),  and that \(r_{k,\tau} = z^{f}\) and \(r_{m,\tau} = z^{2f}\) for \(\tau\) in \(e_{1}ve_{2}\).   Therefore \(r_{m,\tau} = r_{k,\tau}\cdot r_{k,\tau}\) for all \(\tau\in \Sigma\).

\medskip

 \noindent\textsc{Case 2:} Write \(v\) for a vertex of valency 4 marked with \(\chi_{m} = \chi_{k}\otimes \chi_{l}\) from Lemma~\ref{lemma:val4}.  Mimicing the proof of \textsc{Case 1} above gives \(r_{k,\tau} = y^{b}\), \(r_{l,\tau} = y^{c}\), \(r_{m,\tau} = y^{b+c}\) and hence \(r_{m,\tau} = r_{k,\tau}\cdot r_{l,\tau}\) for \(\tau\) in the region \(e_{1}ve_{3}\) of Figure~\ref{fig:Fval4}.  By permuting \(x\) and \(y\),  the same argument shows that \(r_{k,\tau} = x^{d}\), \(r_{l,\tau} = x^{a}\),  \(r_{m,\tau} = x^{a+d}\) and hence \(r_{m,\tau} = r_{k,\tau}\cdot r_{l,\tau}\) for \(\tau\) in \(e_{2}ve_{3}\).  Again,  note in passing that if \(\tau\) has \(v\) as a vertex then \(r_{m,\tau}\) (which is either \(y^{b+c}\) or \(x^{a+d}\)) lies in the socle of \(\Gamma_{\tau}\).

\begin{figure}[!ht]
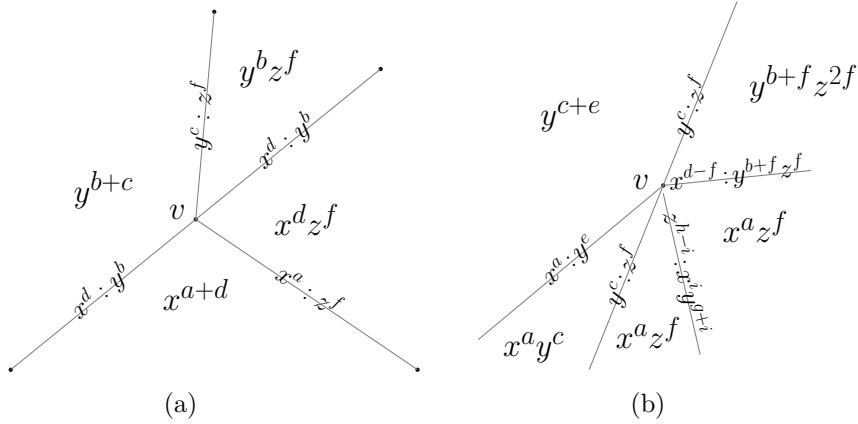

 \centering
 \mbox{\subfigure[]{\input{Fval4m.pspic}\label{fig:Fval4m}} \hspace{1.8cm} \subfigure[]{\input{Fval5m.pspic}\label{fig:Fval5m}}}
 \caption{Generators $r_{m,T}$ for $T$ with vertex $v$ in (a) \textsc{Case 2}; (b) \textsc{Case 3}}
 \end{figure}

 As for \(e_{1}ve_{2}\),  let \(T\) and \(T'\) denote the basic triangles in \(e_{1}ve_{2}\) having \(v\) as a vertex.  Using the coordinates \(\xi = x^{a}/z^{f}, \eta = y^{b}/x^{d}, \zeta = z^{f+1}/y^{b-1}x^{a-d-1}\) on \(U_{T}\) and \(\xi' = x^{d}/y^{b}\),  \(\eta' = y^{c}/z^{f}\), \(\zeta' = z^{f+1}/x^{d-1}y^{c-b-1}\) on \(U_{T'}\),  we calculate the \(A\)-graphs \(\Gamma_{T}\) and \(\Gamma_{T'}\).  It is immediate that \(x^{d}z^{f}\in \Gamma_{T}\) and \(y^{b}z^{f}\in \Gamma_{T'}\) (again,  these monomials both lie in the socle of the corresponding \(A\)-graph).  Since \(\wt(z^{f}) = \chi_{l}\) and \(\wt(x^{d}) = \wt(y^{b}) = \chi_{k}\) it follows that \(r_{m,T} = x^{d}z^{f}\) and \(r_{m,T'} = y^{b}z^{f}\) (see Figure~\ref{fig:Fval4m}).  Lemma~\ref{lemma:facts}\three\ reveals that \(z^{f}\) divides \(r_{m,\tau}\) for every \(\tau\) in \(e_{1}ve_{2}\).  Both \(r_{m,\tau}/z^{f}\) and \(z^{f}\) are of the form \(r_{i,\tau}\) for some \(\chi_{i}\) because \(A\)-graphs are cyclic \(\C[x,y,z]\)-modules with generator 1.  Now,  \(\wt(z^{f}) = \chi_{l}\) and \(\chi_{m} = \chi_{k}\otimes \chi_{l}\),  hence \(\wt(r_{m,\tau}/z^{f}) = \chi_{k}\).  As a result \(r_{l,\tau} = z^{f}\) and \(r_{m,\tau}/z^{f} = r_{k,\tau}\) for \(\tau\) in \(e_{1}ve_{2}\).  This proves that \(r_{m,\tau} = r_{k,\tau}\cdot r_{l,\tau}\) for every triangle \(\tau\) in \(e_{1}ve_{2}\) which completes \textsc{Case 2}.

 \medskip
 \noindent\textsc{Case 3:} Write \(v\) for a vertex of valency 5 or 6 (excluding three straight lines meeting at a point) marked with \(\chi_{m} = \chi_{k}\otimes \chi_{l}\) from Lemma~\ref{lemma:val5or6}.  Mimicing the proof of \textsc{Case 1} above gives \(r_{k,\tau} = y^{c}\),  \(r_{l,\tau} = y^{e}\),  \(r_{m,\tau} = y^{c+e}\) and hence \(r_{m,\tau} = r_{k,\tau}\cdot r_{l,\tau}\) for \(\tau\) in the region \(e_{1}ve_{3}\) of Figure~\ref{fig:Fval5}.  Again,  the monomial \(r_{m,\tau} = y^{c+e}\) lies in the socle of the \(A\)-graph \(\Gamma_{\tau}\) for the triangle \(\tau\) with vertex \(v\).

 To prove the result for \(\tau\) lying outside \(e_{1}ve_{3}\), mimic the proof of \textsc{Case 2} twice.  That is,  first compute the \(A\)-graphs of all triangles \(T\) having \(v\) as a vertex by calculating coordinates on the affine pieces \(U_{T}\).  There exist two such adjacent triangles \(T,T'\) for which \(r_{m,T} = y^{b+f}z^{2f}\) and \(r_{m,T'} = x^{a}z^{f}\) (see Figure~\ref{fig:Fval5m}; these monomials both lie in the socle of the corresponding \(A\)-graph).  Lemma~\ref{lemma:facts}\three\ reveals that \(z^{f}\) divides \(r_{m,\tau}\) for every \(\tau\) whose interior intersects the region \(e_{1}ve_{2}\).   As in \textsc{Case 2},  it follows that \(r_{k,\tau} = z^{f}\),  \(r_{m,\tau}/z^{f} = r_{l,\tau}\) and hence that \(r_{m,\tau} = r_{k,\tau}\cdot r_{l,\tau}\) for \(\tau\) with interior intersecting \(e_{1}ve_{2}\).  Repeat this argument beginning with the adjacent triangles \(T'',T'''\) from Figure~\ref{fig:Fval5m} where \(r_{m,T''} = x^{a}z^{f}\) and \(r_{m,T'''} = x^{a}y^{c}\) to see that \(r_{l,\tau}= x^{a}\),  \(r_{m,\tau}/x^{a} = r_{k,\tau}\) and hence that \(r_{m,\tau} = r_{k,\tau}\cdot r_{l,\tau}\) for \(\tau\) with interior intersecting \(e_{2}ve_{3}\).  The interior of every triangle lying outside \(e_{1}ve_{3}\) intersects either \(e_{1}ve_{2}\) or \(e_{2}ve_{3}\) (possibly both,  in which case \(r_{m,\tau} = x^{a} z^{f}\)),  hence \(r_{m,\tau} = r_{k,\tau}\cdot r_{l,\tau}\) for all \(\tau\) lying outside \(e_{1}ve_{3}\).  This completes the proof of \textsc{Case 3} for a vertex of valency 5 when the regular triangle \(R\) from the proof of Lemma~\ref{lemma:val5or6} is from \(e_{1}\) as in Figure~\ref{fig:Fval5}. The same argument proves \textsc{Case 3} when \(R\) is either a corner triangle from \(e_{2}\) or the meeting of champions.  The case where the vertex has valency 6 is almost identical.

 \medskip

\begin{figure}[!ht]
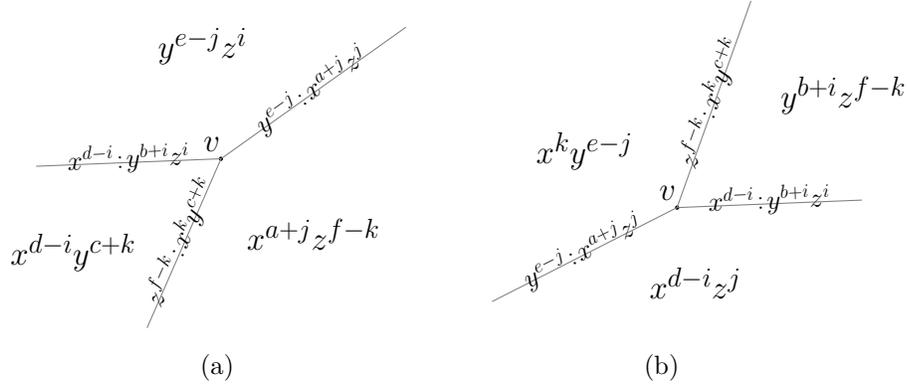

 \centering
 \mbox{\subfigure[]{\input{Fval6l.pspic}\label{fig:Fval6l}} \hspace{1.5cm} \subfigure[]{\input{Fval6m.pspic}\label{fig:Fval6m}}}
 \caption{For $T$ with vertex $v$, the generators of (a) $\mathcal{R}_{l}$; (b) $\mathcal{R}_{m}$}
  \label{fig:Fval6}
 \end{figure}

 \noindent\textsc{Case 4:} The intersection point \(v\) of three straight lines is marked with characters \(\chi_{l}\) and \(\chi_{m}\) from Lemma~\ref{lemma:val6}.  If \(v\) lies inside a corner triangle whose tesselating lines are cut out by the ratios (\ref{eqn:prop3.2a}) then the characters \(\chi_{i}\),  \(\chi_{j}\) and \(\chi_{k}\) mark the lines through \(v\) cut out by \(x^{d-i}\!:\!y^{b+i}z^{i}\),  \(y^{e-j}\!:\!z^{j}x^{a+j}\) and \(z^{f-k}\! : \! x^{k}y^{c+k}\) respectively.  First,  calculate the \(A\)-graphs \(\Gamma_{T}\) of the six triangles \(T\) with \(v\) as a vertex,  and hence write down the generators \(r_{l,T}\) and \(r_{m,T}\) of \(\mathcal{R}_{l}\) and \(\mathcal{R}_{m}\) on the open sets \(U_{T}\) as shown in Figure~\ref{fig:Fval6} (as in the previous three cases,  the monomials shown on Figure~\ref{fig:Fval6} lie in the socle of the corresponding \(A\)-graphs).  The lines in the figure should pass straight through \(v\) leaving six triangles;  for instance,  \(r_{l,T} = y^{e-j}z^{i}\) for a pair of triangles \(T\) with vertex \(v\) (the six monomials written on Figure~\ref{fig:Fval6} first appeared in (\ref{eqn:mapstoP2})).
 
 We claim that every \(\tau\in \Sigma\) whose interior intersects \(e_{1}ve_{2}\) satisfies
 \begin{equation}
 r_{k,\tau} = z^{f-k}, \quad  r_{l,\tau} = z^{i}\cdot r_{j,\tau} \quad \mbox{ and } \quad r_{m,\tau} = z^{j}\cdot r_{i,\tau}. 
 \label{eqn:zf-k}
 \end{equation}
 To see this,  consider the adjacent triangles \(T,T'\) with vertex \(v\) from Figure~\ref{fig:Fval6l} such that \(r_{l,T} = y^{e-j}z^{i}\) and \(r_{l,T'} = x^{a+j}z^{f-k}\).   It follows that \(r_{k,T'} = z^{f-k}\) and hence,  by Lemma~\ref{lemma:facts}\two,  \(r_{k,\tau} = z^{f-k}\) for every \(\tau\) whose interior intersects \(e_{1}ve_{2}\). Moreover,  since \(i < f-k\),  Lemma~\ref{lemma:facts}\three\ reveals that \(z^{i}\) divides \(r_{l,\tau}\) for every \(\tau\) whose interior intersects \(e_{1}ve_{2}\),  hence \(r_{l,\tau}/z^{i}\in \Gamma_{\tau}\).  Observe from (\ref{eqn:mapstoP2}) that \(y^{e-j}z^{i}\) lies in the \(\chi_{l}\)-character space.  This means that \(\wt(r_{l,\tau}/z^{i}) = \wt(y^{e-j}) = \chi_{j}\) and hence \(r_{l,\tau}/z^{i} = r_{j,\tau}\) which proves the second part of (\ref{eqn:zf-k}).  As for the third part,  consider the adjacent triangles \(T',T''\) with vertex \(v\) shown in Figure~\ref{fig:Fval6m} such that \(r_{m,T'} = y^{b+i}z^{f-k}\) and \(r_{m,T''} = x^{d-i}z^{j}\).   Since \(j < f-k\),  Lemma~\ref{lemma:facts}\three\ reveals that \(z^{j}\) divides \(r_{m,\tau}\) for every \(\tau\) with a vertex inside \(e_{1}ve_{2}\),  hence \(r_{m,\tau}/z^{j}\in \Gamma_{\tau}\).  Again,  (\ref{eqn:mapstoP2}) shows that \(x^{d-i}z^{j}\) lies in the \(\chi_{m}\)-character space which means that \(\wt(r_{m,\tau}/z^{j}) = \wt(x^{d-i}) = \chi_{i}\).  Thus \(r_{l,\tau}/z^{i} = r_{i,\tau}\) which completes the proof of (\ref{eqn:zf-k}).  The equations (\ref{eqn:prop3.1a}) show that \(z^{f-k} = z^{i}\cdot z^{j}\),  so (\ref{eqn:zf-k}) gives
 \begin{equation}
 \label{eqn:val6gens}
 r_{l,\tau}\cdot r_{m,\tau} = r_{i,\tau}\cdot r_{j,\tau}\cdot r_{k,\tau}
 \end{equation}
 for all triangles \(\tau\) whose interiors intersect \(e_{1}ve_{2}\).  The next step is to prove that every \(\tau\in \Sigma\) whose interior intersects \(e_{1}ve_{3}\) satisfies
 \begin{equation}
  r_{i,\tau} = x^{d-i},\quad  r_{l,\tau} = x^{a+j}\cdot r_{k,\tau} \quad \mbox{and} \quad r_{m,\tau} = x^{k}\cdot r_{j,\tau}.
  \label{eqn:ye-j}
 \end{equation}
 The proof is the same as that for (\ref{eqn:zf-k}).  In this case, (\ref{eqn:prop3.1a}) gives \(x^{d-i} = x^{a+j}\cdot x^{k}\) which establishes (\ref{eqn:val6gens}) for \(\tau\) whose interiors intersect \(e_{2}ve_{3}\).  Finally,
 \begin{equation}
  r_{j,\tau} = y^{e-j}, \quad r_{l,\tau} = y^{c+k}\cdot r_{i,\tau} \quad \mbox{and} \quad r_{m,\tau} = y^{b+i}\cdot r_{k,\tau} \label{eqn:xd-i}
 \end{equation}
 for every \(\tau\in \Sigma\) whose interior intersects \(e_{1}ve_{3}\).  Again,  the proof is similar to (\ref{eqn:zf-k}), and (\ref{eqn:prop3.1a}) gives \(y^{e-j} = y^{c+k}\cdot y^{b+i}\) thereby proving (\ref{eqn:val6gens}) for \(\tau\) whose interiors intersect \(e_{1}ve_{3}\).  Thus (\ref{eqn:val6gens}) holds for all \(\tau \in \Sigma\) and hence proves \textsc{Case 4} when \(v\) lies inside a corner triangle with tesselating lines cut out by (\ref{eqn:prop3.2a}).  Minor changes in indices proves the case where \(v\) lies inside the meeting of champions triangle whose tesselating lines are cut out by (\ref{eqn:prop3.2b}). 
 \qed

 \begin{remark}
 \emph{In the course of the proof we saw that whenever \(\chi_{m}\) marks \(v\),  the corresponding monomial \(r_{m,\tau}\) lies in the socle of \(\Gamma_{\tau}\) for every triangle having \(v\) as a vertex.  For more on this point see Craw and Ishii~\cite[\S7]{Craw:fog}. 
 }
 \end{remark}

 \begin{proposition}
 \label{prop:multiplicative}
 Theorem~\ref{thm:bundles} lists all nontrivial relations between tautological bundles in \(\Pic(Y)\).  In particular, the map \(\chi_{k}\mapsto \mathcal{R}_{k}\) is not multiplicative.
 \end{proposition}
\proof The nontrivial tautological bundles span \(\Pic(Y)\) because the whole collection \(\{\mathcal{R}_{k}\}\) base \(K(Y)\).  For each relation listed in the statement of the theorem,  remove the bundle \(\mathcal{R}_{m}\) from the spanning set (see Remark~\ref{remark:ambiguity}).  Since we choose \(\mathcal{R}_{m}\) each time,  the remaining tautological bundles indexed by characters of type \one\ and \three\ in Corollary~\ref{coro:chars} still span \(\Pic(Y)\).  There are \(\vert A\vert -1\) nontrivial bundles and we've just removed \(b_{4}(Y)\) of them,  one for each interior vertex in \(\Sigma\),  so the set spanning \(\Pic(Y)\) consists of \(\vert A\vert -1 - b_{4}(Y)\) bundles.  The McKay correspondence of Ito and Reid~\cite{Ito:mc3} gives \(e(Y) = \vert A\vert\),  so the set spanning \(\Pic(Y)\) consists of \(e(Y) -1 - b_{4}(Y)= b_{2}(Y) =\rank\Pic(Y)\) elements.  Thus the bundles indexed by characters of type \one\ and \three\ are independent so there can be no more nontrivial relations. \qed

 \begin{remark}
 \label{remark:ambiguity}
 \emph{There are three maps \(\dP\to \mathbb{P}^{1}\) given by restriction of the bundles \(\mathcal{R}_{i}\), \(\mathcal{R}_{j}\), \(\mathcal{R}_{k}\),  and two maps \(\dP\to \mathbb{P}^{2}\) given by restriction of \(\mathcal{R}_{l}\) and \(\mathcal{R}_{m}\).  All five maps span \(\Pic(\dP)\) and the relation of Theorem~\ref{thm:bundles} part \emph{(4)} holds.  We break the symmetry in this relation by choosing the restriction of the bundles \(\mathcal{R}_{i}\), \(\mathcal{R}_{j}\), \(\mathcal{R}_{k}\), \(\mathcal{R}_{l}\) as a basis for \(\Pic(\dP)\) while discarding \(\mathcal{R}_{m}\).  However,  we could equally well choose the restriction of \(\mathcal{R}_{m}\) as the fourth basis element in which case we discard \(\mathcal{R}_{l}\).}
 \end{remark}

 \section{The McKay correspondence}
 \label{sec:rrdb}
 The last line in the proof of Proposition~\ref{prop:multiplicative} established that the bundles \(\mathcal{R}_{k}\) indexed by characters of type \one\ and \three\ from Corollary~\ref{coro:chars} base \(\Pic(Y)\).  Applying the first Chern class isomorphism gives:

 \begin{proposition}
 \label{prop:H2}
 The first Chern classes \(c_{1}(\mathcal{R}_{k})\) of bundles indexed by characters of type \one\ and \three\ in Corollary~\ref{coro:chars} form a basis of \(H^{2}(Y,\Z)\).
 \end{proposition}

 Next consider \(H^{4}(Y,\Z)\).  Following Reid~\cite{Reid:mc},  we use the relations from Theorem~\ref{thm:bundles} to construct virtual bundles \(\mathcal{V}_{m}\) on \(Y\) indexed by characters \(\chi_{m}\) of type \two\ which, by construction,  have trivial rank and trivial first Chern class.  As before,  we proceed case by case:

 \medskip

 \noindent\textsc{Case 1:} For each relation \(\mathcal{R}_{m} = \mathcal{R}_{k}\otimes\mathcal{R}_{k}\) arising from the marking of a vertex of valency 3,  define \(\mathcal{V}_{m} := \big{(}\mathcal{R}_{k}\oplus \mathcal{R}_{k}\big{)}\ominus \big{(}\mathcal{R}_{m}\oplus \owe_{Y}\big{)}\).

\medskip

 \noindent  \textsc{Case 2:} For each relation \(\mathcal{R}_{m} = \mathcal{R}_{k}\otimes\mathcal{R}_{l}\) arising from the marking of a vertex of valency 4,  define \(\mathcal{V}_{m} := \big{(}\mathcal{R}_{k}\oplus \mathcal{R}_{l}\big{)}\ominus \big{(}\mathcal{R}_{m}\oplus \owe_{Y}\big{)}\).

 \medskip

 \noindent\textsc{Case 3:} For each relation \(\mathcal{R}_{m} = \mathcal{R}_{k}\otimes\mathcal{R}_{l}\) arising from the marking of a vertex of valency 5 or 6,  define \(\mathcal{V}_{m} := \big{(}\mathcal{R}_{k}\oplus \mathcal{R}_{l}\big{)}\ominus \big{(}\mathcal{R}_{m}\oplus \owe_{Y}\big{)}\).
 \medskip

 \noindent\textsc{Case 4:} For each relation \(\mathcal{R}_{l}\otimes\mathcal{R}_{m} = \mathcal{R}_{i}\otimes\mathcal{R}_{j}\otimes\mathcal{R}_{k}\) arising from the marking of the intersection point of three straight lines,  define the virtual bundle \(\mathcal{V}_{m} := \big{(}\mathcal{R}_{i}\oplus \mathcal{R}_{j}\oplus\mathcal{R}_{k}\big{)}\ominus \big{(}\mathcal{R}_{l}\oplus \mathcal{R}_{m}\oplus\owe_{Y}\big{)}\).

 \begin{lemma} 
 \label{lemma:deg1}
 The bundle \(\mathcal{R}_{k}\) has degree one on each curve in \(Y\) defined by a line in \(\Sigma\) marked with \(\chi_{k}\).
 \end{lemma}
 \proof Consider a line in \(\Sigma\) marked with \(\chi_{k}\) and suppose,  permuting \(x,y,z\) if necessary,  that the corresponding curve \(\mathbb{P}^{1}\subset Y\) is parametrised by the ratio \(z^{f-k}\! : \! x^{k}y^{c+k}\).  Then \(\zeta = z^{f-k}/x^{k}y^{c+k}\) is a coordinate on \(U_{T}\) defined by the triangle \(T\) on one side of the line,  hence \(r_{k,T} = x^{k}y^{c+k}\).  Similarly,  \(\nu = x^{k}y^{c+k}/z^{f-k}\) is a coordinate on \(U_{T'}\) defined by the triangle \(T'\) on the other side of the line,  hence \(r_{k,T'} = z^{f-k}\).  Thus the transition function of \(\mathcal{R}_{k}\) on \(U_{T}\cap U_{T'}\) is determined by \(z^{f-k} = \zeta\cdot x^{k}y^{c+k}\).  Since the curve is cut out by \(\zeta = 0\) in \(U_{T}\) we see that \(\mathcal{R}_{k}\) has degree one on the curve.\qed
 
 \begin{proposition}
 \label{prop:H4}
 The classes \(c_{2}(\mathcal{V}_{m})\) form a basis of \(H^{4}(Y,\Z)\) dual to the basis \([S]\in H_{4}(Y,\Z)\) defined by the compact exceptional surfaces \(S\) of the resolution \(\varphi\colon Y \to X\).
 \end{proposition}
 \proof
 The \(\C^{*}\)-action \( (x,y,z)\rightarrow (\lambda x, \lambda y,\lambda z)\) defines a retraction of \(Y\) onto the compactly supported exceptional locus of \(\varphi\),  so the homology classes of the compact exceptional surfaces form an integral basis of \(H_{4}(Y,\Z)\).  Write \(S_{n}\) for the exceptional surface corresponding to the vertex \(v := v_{n}\) in \(\Sigma\) marked with the character \(\chi_{n}\) according to \S\ref{sec:df}. We prove case-by-case that
 \begin{equation}
 \label{eqn:dual}
 \int_{S_{n}} c_{2}(\mathcal{V}_{m}) = \delta_{mn}.
 \end{equation}

 \noindent\textsc{Case 1:} Recall from \textsc{Case 1} of Theorem~\ref{thm:bundles} that \(r_{k,\tau}\)  is \(y^{b}\) (resp.\ \(x^{a}\) or \(z^{f}\)) for \(\tau\) in \(e_{1}ve_{3}\) (resp.\ \(e_{2}ve_{3}\) or \(e_{1}ve_{2}\)),  so \(\mathcal{R}_{k}\) has degree zero on all lines not marked with \(\chi_{k}\).   Also,  \(\mathcal{R}_{k}\) has degree one on the curves defined by a line marked with \(\chi_{k}\) by Lemma~\ref{lemma:deg1}.  Now,  \(\chi_{m}\) marks the vertex \(v_{m}\) of valency 3 corresponding to \(S_{m} = \mathbb{P}^{2}\) and \(\mathcal{R}_{k}\vert_{S_{m}} = \owe_{\mathbb{P}^{2}}(1)\), so
 \[
 \int_{S_{m}} c_{2}(\mathcal{V}_{m}) = \int_{S_{m}} c_{1}\big{(}\mathcal{R}_{k}\vert_{S_{m}}\big{)}^{2} = \owe_{\mathbb{P}^{2}}(1)^{2} = 1,
 \]
 as required.  Next,  consider a vertex \(v_{n}\neq v_{m}\) in \(\Sigma\).  If \(v_{n}\) lies on a line from \(v_{m}\) to some \(e_{j}\) then \(S_{n}\) is a (possibly once or twice blown up) scroll \(\mathbb{F}_{r}\).  The bundle \(\mathcal{R}_{l}\) has degree one on the classes\footnote{We adopt the following notation:  let \(F\), \(M\) and \(D\) denote the classes on a surface scroll \(\mathbb{F}_{r}\) with selfintersection 0, \(r\) and \(-r\) respectively;  we use the same notation for the strict transforms of these classes in a once or twice blown up scroll.} \(M\) and \(D\),  and degree zero on \(F\) (and on each \(-1\)-curve \(E\) if the scroll has been blown up).  Thus \(\mathcal{R}_{k}\vert_{S_{n}} = \owe_{S_{n}}(F) \implies c_{1}(\mathcal{R}_{k}\vert_{S_{n}})^{2} = F^{2} = 0\).  Otherwise \(v_{n}\) lies inside one of the regions \(e_{i}ve_{j}\) in which case \(\mathcal{R}_{k}\vert_{S_{n}} = \owe_{S_{n}}\implies c_{1}(\mathcal{R}_{k}\vert_{S_{n}})^{2} = 0\).  Hence
 \[
 \int_{S_{n}} c_{2}(\mathcal{V}_{m}) = \int_{S_{n}} c_{1}\big{(}\mathcal{R}_{l}\vert_{S_{n}}\big{)}^{2} = 0
 \]
 for any \(v_{n}\neq v_{m}\).  This establishes relation (\ref{eqn:dual}) for \textsc{Case 1}.

 \medskip

 \noindent\textsc{Case 2:} In the notation of \textsc{Case 2} from \S\ref{sec:df},  \(\chi_{k}\) marks the straight line through the vertex \(v_{m}\) of valency 4 so \(\mathcal{R}_{k}\) has degree one on the classes \(M\) and \(D\) on the surface \(S_{m} = \mathbb{F}_{r}\) corresponding to \(v_{m}\).  It follows that \(\mathcal{R}_{k}\vert_{S_{m}} = \owe_{S_{m}}(F)\). Also,  \(\mathcal{R}_{l}\) has degree one on \(F\) because \(\chi_{l}\) marks the other two lines meeting at \(v_{m}\),  so \(\mathcal{R}_{l}\vert_{S_{m}} = \owe_{S_{m}}(M + c\cdot F)\),  for some \(c\in \Z\).  Thus
 \[
 \int_{S_{m}} c_{2}(\mathcal{V}_{m}) = \int_{S_{m}} c_{1}\big{(}\mathcal{R}_{k}\vert_{S_{m}}\big{)}\cdot c_{1}\big{(}\mathcal{R}_{l}\vert_{S_{m}}\big{)} = F\cdot(M + c F) = 1.
 \]
 Next,  consider a vertex \(v_{n}\neq v_{m}\) in \(\Sigma\).  From Figure~\ref{fig:Fval4} we see that \(\mathcal{R}_{k}\) (and \(\mathcal{R}_{l}\)) has degree one (resp.\ degree \(d\geq 1\)) on the line \(v_{m}\) to \(e_{3}\),  and degree zero (resp.\ degree 1) on the lines \(v_{m}\) to \(e_{1}\) and \(v_{m}\) to \(e_{2}\).  If \(v_{n}\) lies on the line \(v_{m}\) to \(e_{3}\) then \(S_{n}\) is a scroll \(\mathbb{F}_{r}\) (possibly blown up in one or two points) and,  as above,  we have \(\mathcal{R}_{k}\vert_{S_{n}} = \owe_{S_{n}}(F)\) and \(\mathcal{R}_{l}\vert_{S_{n}} = \owe_{S_{n}}(d F)\) for some \(d\in \Z\).  Thus
 \[
 \int_{S_{n}} c_{2}(\mathcal{V}_{m}) = \int_{S_{n}} c_{1}\big{(}\mathcal{R}_{k}\vert_{S_{n}}\big{)}\cdot c_{1}\big{(}\mathcal{R}_{l}\vert_{S_{n}}\big{)} = F\cdot(d F) = 0.
 \]
If \(v_{n}\neq v_{m}\) lies on the line \(v_{m}\) to \(e_{1}\) or \(v_{m}\) to \(e_{2}\) then \(\mathcal{R}_{k}\vert_{S_{n}} = \owe_{S_{n}}\),  and if \(v_{n}\neq v_{m}\) does not lie on a line from \(v_{m}\) to some \(e_{j}\) then \(\mathcal{R}_{l}\vert_{S_{n}} = \owe_{S_{n}}\).  In either case the Chern class calculation is zero so the relation (\ref{eqn:dual}) holds.

 \medskip

 \noindent\textsc{Case 3:} Almost identical to \textsc{Case 2} so we leave it as an exercise.

 \medskip

 \noindent\textsc{Case 4:} In the notation of \textsc{Case 4} from \S\ref{sec:df},  write \(v_{m}\) for the vertex marked with \(\chi_{l}\) and \(\chi_{m}\) defining a surface \(S_{m}:= \dP\).  The divisor class group is
\[\Div(S_{m}) = \left\langle D_{1}, D_{2}, D_{3}, C_{1}, C_{2} \bigm{|} D_{1}+ D_{2}+ D_{3} = C_{1}+ C_{2}\right\rangle,\]
where \(\owe_{S_{m}}(C_{\alpha})\) and \(\owe_{S_{m}}(D_{\beta})\) define morphisms from \(S_{m}\) to \(\mathbb{P}^{2}\) and \(\mathbb{P}^{1}\) respectively.  The characters \(\chi_{i}\), \(\chi_{j}\) and \(\chi_{k}\) mark the straight lines passing through \(v_{m}\) cut out by the ratios (\ref{eqn:prop3.2a}) or (\ref{eqn:prop3.2b}),  and it follows that \(\mathcal{R}_{i}\vert_{S_{m}} = \owe_{S_{m}}(D_{1})\),  \(\mathcal{R}_{j}\vert_{S_{m}} = \owe_{S_{m}}(D_{2})\) and \(\mathcal{R}_{k}\vert_{S_{m}} = \owe_{S_{m}}(D_{3})\).  Thus
 \[
 \int_{S_{m}} c_{2}(\mathcal{R}_{i}\oplus \mathcal{R}_{j}\oplus\mathcal{R}_{k}) = \sum_{\alpha<\beta} D_{\alpha}\cdot D_{\beta} = 3,
 \]
Also,  from the construction of the characters \(\chi_{l}\) and \(\chi_{m}\) in Lemma~\ref{lemma:val6},  we have \(\mathcal{R}_{l}\vert_{S_{m}} = \owe_{S_{m}}(C_{1})\) and \(\mathcal{R}_{m}\vert_{S_{m}} = \owe_{S_{m}}(C_{2})\), so
 \[\int_{S_{m}}c_{2}\big{(}\mathcal{R}_{l}\oplus \mathcal{R}_{m}\big{)} = \int_{S_{m}}c_{1}(\mathcal{R}_{l})\cdot c_{1}(\mathcal{R}_{m}) = C_{1}\cdot C_{2} = 2.\]
The difference of these two integrals establishes (\ref{eqn:dual}) when \(m=n\).  Next,  consider a vertex \(v_{n}\neq v_{m}\).  When \(v_{n}\) lies in the region \(e_{1}ve_{2}\),  the proof of (\ref{eqn:zf-k}) shows that \(r_{k,\tau} = z^{f-k}\) on the open sets \(U_{\tau}\subset S_{n}\) defined by triangles \(\tau\) with \(v_{n}\) as a vertex.  Then \(\mathcal{R}_{k}\vert_{S_{n}} = \owe_{S_{n}}\),  so \(c_{1}(\mathcal{R}_{k}\vert_{S_{n}}) = 0\).  It follows from (\ref{eqn:zf-k}) that \(c_{1}(\mathcal{R}_{l}\vert_{S_{n}}) = c_{1}(\mathcal{R}_{j}\vert_{S_{n}})\) and \(c_{1}(\mathcal{R}_{m}\vert_{S_{n}}) = c_{1}(\mathcal{R}_{i}\vert_{S_{n}})\).  As a result
 \begin{eqnarray*}
 \int_{S_{n}} c_{2}(\mathcal{R}_{i}\oplus\mathcal{R}_{j}\oplus\mathcal{R}_{k}) & = & c_{1}(\mathcal{R}_{i}\vert_{S_{n}})\cdot c_{1}(\mathcal{R}_{j}\vert_{S_{n}}) \\
 & = & c_{1}(\mathcal{R}_{l}\vert_{S_{n}})\cdot c_{1}(\mathcal{R}_{m}\vert_{S_{n}}) = \int_{S_{n}} c_{2}(\mathcal{R}_{l}\oplus\mathcal{R}_{m}),
 \end{eqnarray*}
 so (\ref{eqn:dual}) holds for \(n\neq m\).  This completes the proof of Proposition~\ref{prop:H4}. 
 \qed

 \medskip
 
 \noindent\emph{Proof of Theorem~\ref{thm:McKay}.\ } The basis \(c_{2}(\mathcal{V}_{m})\) of \(H^{4}(Y,\Z)\) is indexed by characters of type \two,   the basis \(c_{1}(\mathcal{R}_{k})\) of \(H^{2}(Y,\Z)\) is indexed by characters of type \one\ and \three,  the trivial bundle \(\mathcal{R}_{0} = \owe_{Y}\) generates \(H^{0}(Y,\Z)\).
 \qed

 \begin{remark}
 \emph{If \(\chi_{l}\) and \(\chi_{m}\) mark the same vertex then there is a choice as to whether \(c_{1}(\mathcal{R}_{l})\) or \(c_{1}(\mathcal{R}_{m})\) is a basis element of \(H^{2}(Y,\Z)\),  and as to whether we label the virtual bundle in \textsc{Case 4} as \(\mathcal{V}_{m}\) or \(\mathcal{V}_{l}\).  In particular,  when there is a del Pezzo surface \(\dP\subset Y\), there is no canonical answer to the question `Which characters of the group correspond to elements of \(H^{4}(Y,\Z)\) and which to \(H^{2}(Y,\Z)\)?'.
 }
 \end{remark}

 \newpage

 \bibliography{research}

 \medskip

 \noindent \textsc{Department of Mathematics,  University of Utah,\\
 155 South 1400 East, Salt Lake City,UT 84112, USA. \\
 E-mail:} \texttt{craw@math.utah.edu}
 
 \end{document}